\documentclass{amsart}

\title[Combinatorics of KP solitons from the real Grassmannian]
{Combinatorics of KP solitons from the real Grassmannian}

\author{Yuji Kodama and Lauren Williams} 
\date{\today}
\thanks{The first author was partially
supported by NSF grants DMS-0806219 and DMS-1108813.  The second author was 
partially supported by an NSF CAREER award and an 
Alfred Sloan Fellowship.}
\address{Department of Mathematics, Ohio State University,
Columbus, OH 43210}
\email{kodama@math.ohio-state.edu}
\address{Department of Mathematics, University of California,
Berkeley, CA 94720-3840}
\email{williams@math.berkeley.edu}
\subjclass[2000]{}

\usepackage{rotating}
\newcommand{\rotxc}[1]{\begin{sideways}#1\end{sideways}}
\newcommand{\invert}[1]{\rotxc{\rotxc{#1}}}

\def\Le{\hbox{\invert{$\Gamma$}}}
\countdef\x=23
\countdef\y=24
\countdef\z=25
\countdef\t=26

\def\tbox(#1,#2)#3{
\x=#1 \y=#2 
\multiply\x by 12 
\multiply\y by 12 
\z=\x \t=\y
\advance\z by 12 
\advance\t by 12 
\psline(\x,\y)(\x,\t)(\z,\t)(\z,\y)(\x,\y)
\advance\x by 6
\advance\y by 6 
\rput(\x,\y){{\bf #3}}}

\usepackage{amssymb}
\usepackage{graphicx}
\usepackage{epstopdf}
\usepackage{amscd}
\usepackage{graphics}
\def\proof{\par{\it Proof}. \ignorespaces}
\def\endproof{{\ \vbox{\hrule\hbox{%
     \vrule height1.3ex\hskip0.8ex\vrule}\hrule }}\par}

\theoremstyle{definition}

\theoremstyle{remark}

\numberwithin{equation}{section}

\let\trueint=\int
\let\truesum=\sum
\def\int{\mathop{\textstyle\trueint}\limits}
\def\sum{\mathop{\textstyle\truesum}\limits}
\let\<=\langle
\let\>=\rangle

\def\F{{\mathbb{F}}}

\def\S{{\mathcal S}}
\def\J{{\mathcal J}}

\def\SL{\mathrm{SL}}
\def\Sym{{\mathfrak S}}

\def\v{\mathbf{v}}
\def\w{\mathbf{w}}
\def\wstn{\raisebox{0.12cm}{\hskip0.14cm\circle{4}\hskip-0.15cm}}
\def\bstn{\raisebox{0.12cm}{\hskip0.14cm\circle*{4}\hskip-0.15cm}}

\usepackage{amsmath, amsthm, amssymb, amsbsy,amsfonts,amscd}
\usepackage[mathscr]{eucal}
\usepackage{epsfig}
\usepackage{young}
\newtheorem{theorem}{Theorem}[section]

\newtheorem{definition}[theorem]{Definition}
\newtheorem{proposition}[theorem]{Proposition}
\newtheorem{lemma}[theorem]{Lemma}

\newtheorem{example}[theorem]{Example}
\newtheorem{corollary}[theorem]{Corollary}

\newtheorem{remark}[theorem]{Remark}
\newtheorem{algorithm}[theorem]{Algorithm}
\usepackage{amsfonts}
\usepackage{xy}
\usepackage{amssymb}
\usepackage{amsmath}
\newcommand{\inv}{^{-1}}
\newcommand{\To}{\longrightarrow}

\newcommand{\R}{\mathbb R}

\newcommand{\B}{\mathcal{B}}
\newcommand{\Grkn}{(Gr_{k,n})_{\geq 0}}

\newcommand{\Q}{\mathcal{Q}}

\DeclareMathOperator{\Dec}{Dec}

\DeclareMathOperator{\CC}{\mathcal C}

\DeclareMathOperator{\In}{in}
\DeclareMathOperator{\Out}{out}

\DeclareMathOperator{\M}{\mathcal M}
\DeclareMathOperator{\I}{\mathcal I}

\newcommand{\thmrefer}[1]{\renewcommand\thetheorem
  {\protect\ref{#1}}\addtocounter{theorem}{-1}}

\xyoption{all}
\CompileMatrices

\begin{document}

\begin{abstract}
Given a point $A$ in the real Grassmannian, it is well-known that
one can construct a soliton solution $u_A(x,y,t)$ to the KP equation. 
The \emph{contour plot} of such a solution provides a tropical 
approximation to the solution when the variables $x$, $y$, and $t$ 
are considered 
on a large scale and the time $t$ is fixed.
In this paper we 
give an overview of our work on the combinatorics of 
such contour plots.  Using the 
positroid stratification and the Deodhar decomposition
of the Grassmannian (and in particular the combinatorics
of \emph{Go-diagrams}), we completely describe 
the asymptotics of these contour plots when 
$y$ or $t$ go to $\pm \infty$. 
Other highlights include: a surprising connection
with total positivity and cluster algebras; results on the 
\emph{inverse problem}; and
the characterization of regular soliton solutions -- 
that is, a soliton solution $u_A(x,y,t)$ is regular for all times $t$ 
if and only if $A$ comes from the \emph{totally non-negative part}
of the Grassmannian.
\end{abstract}

\maketitle


\begin{center}
\begin{figure}[h]
\includegraphics[height=.5in]{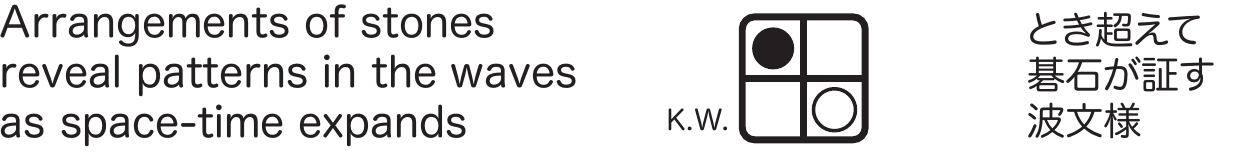}
\end{figure}
\end{center}

\bigskip

\setcounter{tocdepth}{1}
\tableofcontents

\section{Introduction}

The main purpose of this paper is to give an exposition of our recent work
\cite{KW, KW2, KW3},
which found surprising connections between soliton solutions of the KP equation and the combinatorics of the real Grassmannian.
The KP equation is a two-dimensional
nonlinear dispersive wave equation which was proposed by Kadomtsev and Peviashvili
in 1970 to study the stability problem of the soliton solution of the Korteweg-de Vries (KdV)
equation \cite{KP70}.  
The equation has a rich
mathematical structure, and is now considered
to be the prototype of an integrable nonlinear
dispersive wave equation with two spatial dimensions (see for example  \cite{NMPZ84,AC91,D91,MJD00,H04}).
The KP equation can also be used to describe shallow water wave phenomena,
including resonant interactions.

An important breakthrough in the KP theory was made by Sato \cite{Sato}, who
realized that solutions of the KP equation could be written in terms of points on an
infinite-dimensional Grassmannian.  The present paper, which 
gives an overview of 
most of our results of \cite{KW2} and \cite{KW3}, 
deals with a real, finite-dimensional
version of the Sato theory.   In particular, we are interested in 
\emph{soliton solutions}, that is, solutions that
are localized along
certain rays in the $xy$ plane called \emph{line-solitons}.  
Such a 
solution
can be constructed from a point $A$ of the real Grassmannian.  
More specifically,
one can apply
the {\it Wronskian form} \cite{Sato,Sa79,FN,H04}
to $A$ to produce a certain sum of exponentials called 
a \emph{$\tau$-function} $\tau_A(x,y,t)$, and from the
$\tau$-function 
one can construct a solution $u_A(x,y,t)$ to the KP equation.

Recently several authors
have studied the soliton solutions $u_A(x,y,t)$ which come from points $A$ of 
the {\it totally non-negative part of the Grassmannian} $(Gr_{k,n})_{\geq 0}$, that is,
those points of the real Grassmannian $Gr_{k,n}$
whose Pl\"ucker coordinates are all non-negative
\cite{BK,K04, BC06, CK1,CK3, KW, KW2}.  These solutions are 
\emph{regular}, and include
a large variety of soliton solutions which were previously overlooked by those using the
Hirota method of a perturbation expansion \cite{H04}.


A main goal of \cite{KW3} was to understand the soliton solutions
$u_A(x,y,t)$ coming from arbitrary points $A$ of the real Grassmannian, not just 
the totally non-negative part. In general such solutions are no longer regular -- 
they may have singularities along rays in the $xy$ plane -- 
but it is possible, nevertheless, to understand a great deal about the asymptotics of such solutions,
when the absolute value of the spatial variable $y$ goes to infinity, and also when
the absolute value of the time variable $t$ goes to infinity.

Two related decompositions of the real Grassmannian are  useful for understanding
the asymptotics of soliton solutions $u_A(x,y,t)$.
The first is 
Postnikov's \emph{positroid stratification} of the Grassmannian 
\cite{Postnikov}, whose strata are indexed by various combinatorial
objects including decorated permutations and $\Le$-diagrams.  This decomposition
determines the asympototics of soliton solutions when $|y| \gg 0$, and our
results here 
extend work of \cite{BC06, CK1, CK3, KW, KW2}
from the setting of the non-negative part of the Grassmannian to the entire real Grassmannian.
The second decomposition, which we call the \emph{Deodhar decomposition} \cite{KW3},
is the projection to the Grassmannian of 
a decomposition of the flag variety due to Deodhar.  
The Deodhar decomposition refines the positroid stratification,
and its components may be indexed by certain tableaux filled with black
and white stones called \emph{Go-diagrams}, which generalize $\Le$-diagrams.
This decomposition determines the asymptotics of soliton solutions
when $|t| \gg 0$.  More specifically, it allows us to compute the 
\emph{contour plots} at $|t| \gg 0$ of such solitons, which are tropical approximations
to the solution when $x$, $y$, and $t$ are on a large scale \cite{KW3}.

By using our results on the asymptotics of soliton solutions when $t\ll 0$,
one may give a characterization of the regular soliton solutions coming from the real Grassmannian.
More specifically, a soliton solution $u_A(x,y,t)$ coming from a point $A$ of the real Grassmannian
is regular for all times $t$ if and only if $A$ is a point of the totally non-negative part of
the Grassmannian \cite{KW3}.

The regularity theorem above provides an important motivation for studying 
soliton solutions coming from $(Gr_{k,n})_{\geq 0}$.  Indeed, as 
we showed in \cite{KW2}, such soliton solutions
have an even richer combinatorial structure than those coming from $Gr_{k,n}$.
For example, (generic) contour plots coming from the totally
positive part $(Gr_{k,n})_{>0}$ of the Grassmannian give rise to \emph{clusters}
for the cluster algebra associated to the Grassmannian.
And up to a combinatorial equivalence, the contour plots coming from 
$(Gr_{2,n})_{>0}$ are in bijection with triangulations of an $n$-gon.
Finally, if either $A\in (Gr_{k,n})_{>0}$, or 
$A \in (Gr_{k,n})_{\geq 0}$ and $t \ll0$, then one may solve the 
\emph{inverse problem} for $u_A(x,y,t)$: that is, given the 
contour plot $\CC_t(u_A)$ and the time $t$, one may reconstruct the 
element $A \in (Gr_{k,n})_{\geq 0}$.  And therefore one may reconstruct 
the entire evolution of this soliton solution over time \cite{KW2}.

The structure of this paper is as follows.
In Section \ref{sec:Gr} we provide background on the Grassmannian 
and some of its decompositions, including the positroid stratification.
In Section \ref{sec:project} we describe
the Deodhar decomposition of the complete flag variety and its projection
to the Grassmannian, while in Section \ref{Deodhar-combinatorics} we
explain how to index Deodhar components in the Grassmannian by 
\emph{Go-diagrams}.
Subsequent sections provide applications
of the previous results to soliton solutions of the KP equation.
In Section \ref{soliton-background} we explain how to produce
a soliton solution to the KP equation from a point of the real
Grassmannian, and then 
define the \emph{contour plot} associated to a soliton solution
at a fixed time $t$.  In Section \ref{sec:unbounded} we use the positroid stratification to 
describe the unbounded line-solitons in contour plots of soliton solutions
at $y\gg 0$ and $y \ll 0$.
In Section \ref{sec:solitonplabic}
we define the more combinatorial notions of 
\emph{soliton graph} and \emph{generalized plabic graph}.
In Section \ref{sec:t<<0} we use the Deodhar
decomposition to describe contour plots of soliton solutions
for $t\ll 0$.
In Section \ref{sec:regularity} we describe the significance of total
positivity to soliton solutions, by discussing the regularity
problem, as well as the connection to cluster algebras.  Finally
in Section \ref{sec:inverse}, we give results on the inverse problem
for soliton solutions coming from $(Gr_{k,n})_{\geq 0}$.


\section{Background on the Grassmannian and its totally non-negative part}\label{sec:Gr}

The \emph{real Grassmannian} $Gr_{k,n}$ is the space of all
$k$-dimensional subspaces of $\R^n$.  An element of
$Gr_{k,n}$ can be viewed as a full-rank $k\times n$ matrix modulo left
multiplication by nonsingular $k\times k$ matrices.  In other words, two
$k\times n$ matrices represent the same point in $Gr_{k,n}$ if and only if they
can be obtained from each other by row operations.
Let $\binom{[n]}{k}$ be the set of all $k$-element subsets of $[n]:=\{1,\dots,n\}$.
For $I\in \binom{[n]}{k}$, let $\Delta_I(A)$
be the {\it Pl\"ucker coordinate}, that is, the maximal minor of the $k\times n$ matrix $A$ located in the column set $I$.
The map $A\mapsto (\Delta_I(A))$, where $I$ ranges over $\binom{[n]}{k}$,
induces the {\it Pl\"ucker embedding\/} $Gr_{k,n}\hookrightarrow \mathbb{RP}^{\binom{n}{k}-1}$.
The \emph{totally non-negative part of the Grassmannian}
$(Gr_{k,n})_{\geq 0}$ is the subset of $Gr_{k,n}$ such that all 
Pl\"ucker coordinates are non-negative.

We now describe several useful decompositions of the Grassmannian:
the matroid stratification, the Schubert decomposition, and the positroid
stratification.  When one restricts the positroid stratification to 
$(Gr_{k,n})_{\geq 0}$, one gets a cell decomposition
of $(Gr_{k,n})_{\geq 0}$ into \emph{positroid cells}.

\subsection{The matroid stratification of $Gr_{k,n}$}

\begin{definition}\label{def:matroid}
A \emph{matroid} of \emph{rank} $k$ on the set $[n]$ is a nonempty collection
$\M \subset \binom{[n]}{k}$ of $k$-element subsets in $[n]$, called \emph{bases}
of $\M$, that satisfies the \emph{exchange axiom}:\\
For any $I,J \in \M$ and $i \in I$ there exists $j\in J$ such that
$(I \setminus \{i\}) \cup \{j\} \in \M$.
\end{definition}

Given an element $A \in Gr_{k,n}$, there is an associated matroid
$\M_A$ whose bases are the $k$-subsets $I \subset [n]$ such that 
$\Delta_I(A) \neq 0$.

\begin{definition}
Let $\M \subset \binom{[n]}{k}$ be a matroid.
The \emph{matroid stratum}
$S_{\M}$ is defined to be 
$$S_{\M} = \{A \in Gr_{k,n} \ \vert \ \Delta_I(A) \neq 0 \text{ if and only if }
I\in \M \}.$$
This gives a stratification of $Gr_{k,n}$ called the 
\emph{matroid stratification}, or \emph{Gelfand-Serganova stratification}.
The matroids $\M$ with nonempty strata $S_{\M}$ are called \emph{realizable} over
$\R$.
\end{definition}

\subsection{The Schubert decomposition of $Gr_{k,n}$}

Recall that the partitions $\lambda \subset (n-k)^k$
are in bijection with $k$-element subset $I \subset [n]$.  
The boundary of the Young diagram of such a partition 
$\lambda$ forms a lattice path from the upper-right corner to the lower-left
corner of the rectangle $(n-k)^k$.  Let us label the $n$ steps 
in this path by the numbers $1,\dots,n$, and define
$I = I(\lambda)$ as the set of labels on the $k$ vertical steps in 
the path. Conversely, we let $\lambda(I)$ denote the 
partition corresponding to the subset $I$.

\begin{definition}
For each partition $\lambda \subset (n-k)^k$, one can define the 
\emph{Schubert cell} $\Omega_{\lambda}$ by 
$$\Omega_{\lambda} = \{A \in Gr_{k,n} \ \vert \ I(\lambda) \text{ is
the lexicographically minimal base of }\M_A \}.$$
As $\lambda$ ranges over the partitions contained in $(n-k)^k$,
this gives the \emph{Schubert decomposition} 
of the Grassmannian $Gr_{k,n}$, i.e.
\[
Gr_{k,n}=\bigsqcup_{\lambda\subset (n-k)^k}\,\Omega_{\lambda}.
\]
\end{definition}

We now define the \emph{shifted linear order} $<_i$ (for $i\in [n]$) to be the total order on $[n]$
defined by $$i <_i i+1 <_i i+2 <_i \dots <_i n <_i 1 <_i \dots <_i i-1.$$
One can then define \emph{cyclically shifted Schubert cells} as follows.

\begin{definition} \label{def:Schubert}
For each partition $\lambda \subset (n-k)^k$ and $i \in [n]$, the 
\emph{cyclically shifted Schubert cell} $\Omega_{\lambda}^i$ is
$$\Omega_{\lambda}^i = \{A \in Gr_{k,n} \ \vert \ I(\lambda) \text{ is
the lexicographically minimal base of }\M_A \text{ with respect to }<_i \}.$$
\end{definition}


\subsection{The positroid stratification of $Gr_{k,n}$}\label{sec:positroid}

The \emph{positroid stratification} of the real Grassmannian
$Gr_{k,n}$ is obtained by taking the
simultaneous refinement of the $n$ Schubert decompositions with respect to 
the $n$ shifted linear orders $<_i$.  This stratification 
was first considered
by Postnikov \cite{Postnikov}, who showed that
the strata are conveniently described 
in terms of 
\emph{Grassmann necklaces}, as well as \emph{decorated permutations},
(equivalence classes of) \emph{plabic graphs}, 
and \emph{$\Le$-diagrams}.  Postnikov coined the terminology
\emph{positroid} because the intersection of the positroid stratification
with
$(Gr_{k,n})_{\geq 0}$ gives a cell decomposition of 
$(Gr_{k,n})_{\geq 0}$  (whose cells are called 
\emph{positroid cells}).

\begin{definition}\cite[Definition 16.1]{Postnikov}
A \emph{Grassmann necklace} is a sequence
$\I = (I_1,\dots,I_n)$ of subsets $I_r \subset [n]$ such that,
for $i\in [n]$, if $i\in I_i$ then $I_{i+1} = (I_i \setminus \{i\}) \cup \{j\}$,
for some $j\in [n]$; and if $i \notin I_r$ then $I_{i+1} = I_i$.
(Here indices $i$ are taken modulo $n$.)  In particular, we have
$|I_1| = \dots = |I_n|$, which is equal to some $k \in [n]$.  We then say that 
$\I$ is a Grassmann necklace of \emph{type} $(k,n)$.
\end{definition}

\begin{example}\label{ex1}
$(1257, 2357, 3457, 4567, 5678, 6789, 1789, 1289, 1259)$ is
a Grassmann necklace of type $(4,9)$. 
\end{example}

\begin{lemma}\cite[Lemma 16.3]{Postnikov} \label{lem:necklace}
Given $A\in Gr_{k,n}$, 
let $\mathcal I(A) = (I_1,\dots,I_n)$ be the sequence of subsets in 
$[n]$ such that, for $i \in [n]$, $I_i$ is the lexicographically
minimal subset of $\binom{[n]}{k}$  with respect to the shifted linear order
$<_i$ such that $\Delta_{I_i}(A) \neq 0$.
Then $\I(A)$ is a Grassmann necklace of type $(k,n)$.
\end{lemma}

If $A$ is in the matroid stratum $S_{\M}$, 
we also use $\I_{\M}$ to denote the sequence
$(I_1,\dots,I_n)$ defined above.
This leads to the following description of the 
\emph{positroid stratification}
of $Gr_{k,n}$.

\begin{definition}\label{def:pos}
Let $\I = (I_1,\dots,I_n)$ be a Grassmann necklace of type $(k,n)$.
The \emph{positroid stratum}
$S_{\I}$ is defined to be 
\begin{align*}
S_{\I} &= \{A \in Gr_{k,n} \ \vert \ \I(A) = \I \} \\
  &= \bigcap_{i=1}^n ~\Omega_{\lambda(I_i)}^i\,.
\end{align*}
\end{definition}
The second equality follows from Definition \ref{def:pos} and
Definition \ref{def:Schubert}.
Note that  each positroid stratum is an intersection of $n$ 
cyclically shifted Schubert cells.


\begin{definition}\cite[Definition 13.3]{Postnikov}
A \emph{decorated permutation} $\pi^{:} = (\pi, col)$
is a permutation $\pi \in S_n$ together with a coloring
function $col$ from the set of fixed points
$\{i \ \vert \ \pi(i) = i\}$ to $\{1,-1\}$.  So
a decorated permutation is a permutation with fixed points
colored in one of two colors.  A \emph{weak excedance} of
$\pi^{:}$ is a pair $(i,\pi_i)$  such that either
$\pi(i)>i$, or $\pi(i)=i$ and $col(i)=1$. We call
$i$ the \emph{weak excedance position}.
If $\pi(i)>i$ (respectively $\pi(i)<i$) then 
$(i,\pi_i)$ is called an excedance (respectively, nonexcedance).
\end{definition}

\begin{example}\label{ex2}
The decorated permutation (written in one-line notation)
$(6,7,1,2,8,3,9,4,5)$ has no fixed points,
and 
four weak excedances, in positions $1, 2, 5$ and $7$.
\end{example}

\begin{figure}[h]
\centering
\includegraphics[height=1.25in]{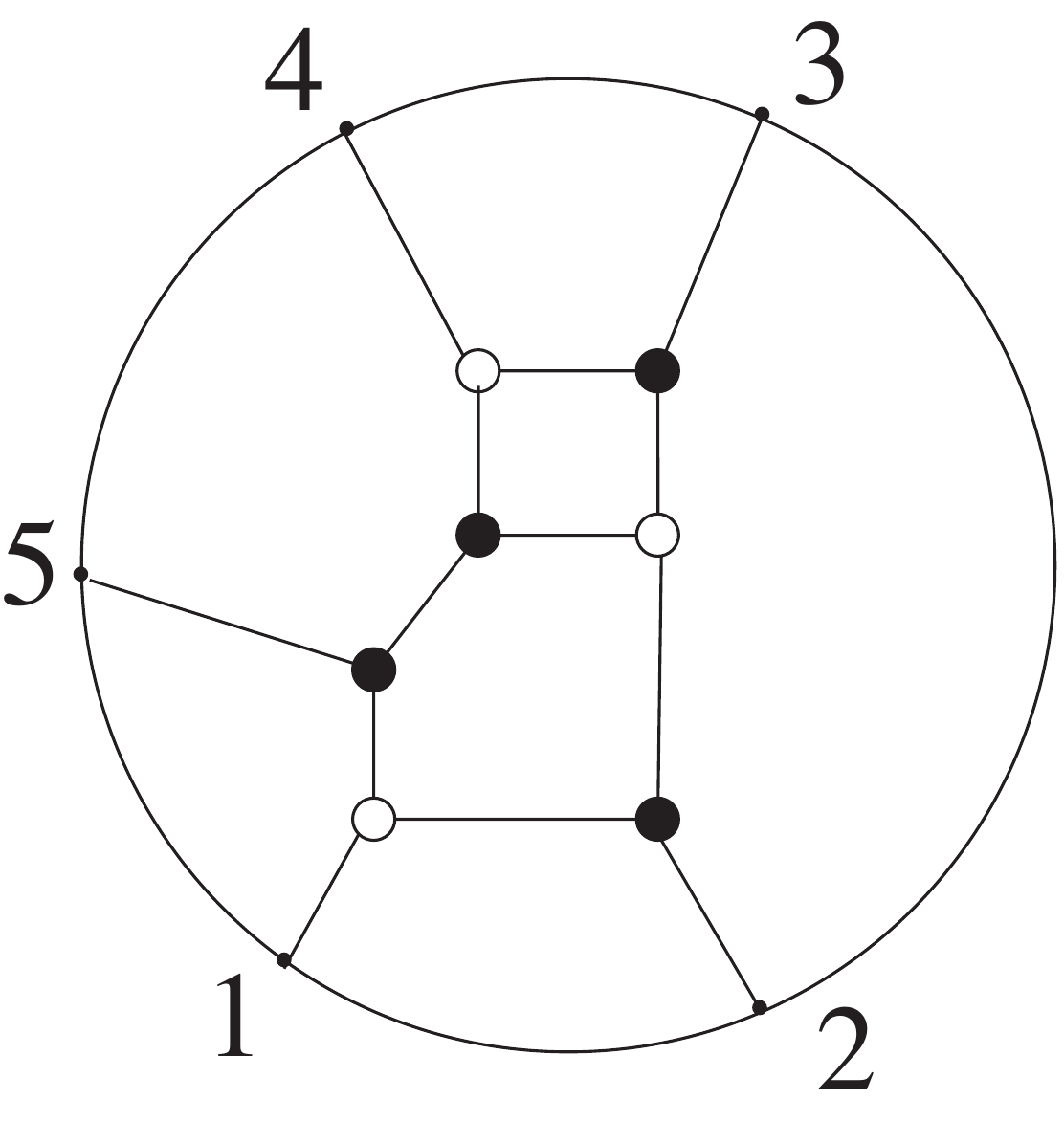}
\hskip 1cm
\includegraphics[height=1.2in]{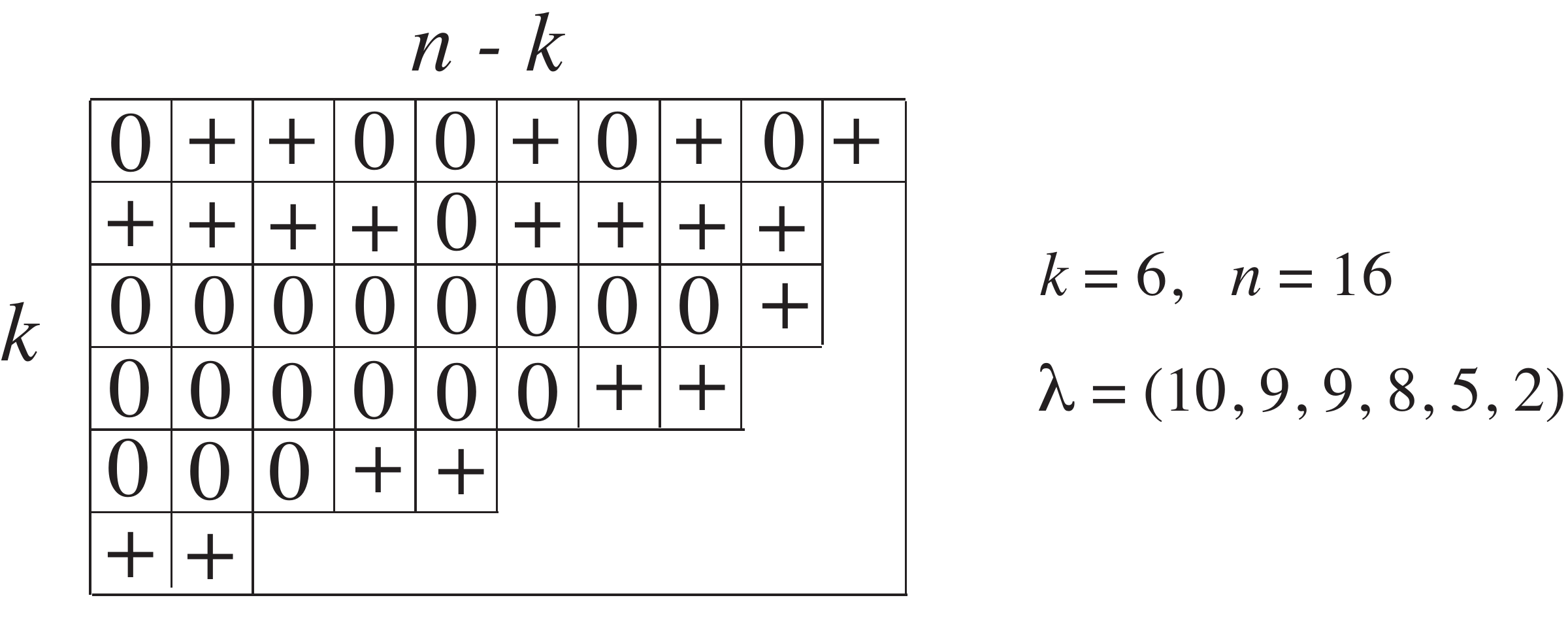}
\caption{A plabic graph and a Le-diagram $L=(\lambda,D)_{k,n}$.
\label{LeDiagram}}
\end{figure}

\begin{definition}\label{def:plabic}
A \emph{plabic graph}
is a planar undirected graph $G$ drawn inside a disk
with $n$ \emph{boundary vertices\/} $1,\dots,n$ placed in counterclockwise
order
around the boundary of the disk, such that each boundary vertex $i$
is incident
to a single edge.\footnote{The convention of \cite{Postnikov}
was to place the boundary vertices in clockwise order.}  Each internal vertex
is colored black or white.  See the left of Figure \ref{LeDiagram} for an example.
\end{definition}

\begin{definition}\cite[Definition 6.1]{Postnikov}\label{def:Le}
Fix $k$, $n$. If $\lambda$ is a partition, let
$Y_{\lambda}$ denote its Young diagram.  A {\it $\Le$-diagram}
$(\lambda, D)_{k,n}$ of type $(k,n)$
is a partition $\lambda$ contained in a $k \times (n-k)$ rectangle
together with a filling $D: Y_{\lambda} \to \{0,+\}$ which has the
{\it $\Le$-property}: 
there is no $0$ which has a $+$ above it and a $+$ to its
left.\footnote{This forbidden pattern is in the shape of a backwards $L$,
and hence is denoted $\Le$ and pronounced ``Le."}  (Here, ``above" means above and in the same column, and
``to its left" means to the left and in the same row.)
See the right of  Figure \ref{LeDiagram} for an example.
\end{definition}

We now review some of the bijections among these objects.

\begin{definition}\cite[Section 16]{Postnikov}\label{necklace-to-perm}
Given a Grassmann necklace $\mathcal I$, define
a decorated permutation $\pi^{:}=\pi^{:}(\mathcal I)$ by requiring that
\begin{enumerate}
\item if $I_{i+1} = (I_i \setminus \{i\}) \cup \{j\}$,
for $j \neq i$, then $\pi^{:}(j)=i$.
\footnote{Actually Postnikov's convention was to set $\pi^:(i)=j$ above,
so the decorated permutation we are associating is the inverse one to his.}
\item if $I_{i+1}=I_i$ and $i \in I_i$ then $\pi(i)=i$ is colored
 with $col(i)=1$.
\item if $I_{i+1}=I_i$ and $i \notin I_i$ then $\pi(i)=i$ is colored
 with $col(i)=-1$.
\end{enumerate}
As before, indices are taken modulo $n$.
\end{definition}
If $\pi^{:}=\pi^{:}(\mathcal I)$, then we also use
the notation
$S_{\pi^{:}}$ to refer to the positroid stratum
$S_{\I}$.  


\begin{example}\label{ex3}
Definition \ref{necklace-to-perm} carries the Grassmann necklace
of Example \ref{ex1} to the decorated permutation of Example \ref{ex2}.
\end{example}

\begin{lemma}\cite[Lemma 16.2]{Postnikov}\label{Postnikov-permutation}
The map $\mathcal I \to \pi^{:}(\mathcal I)$
is a bijection from
Grassmann necklaces $\mathcal I=(I_1,\dots,I_n)$
of size $n$ to
decorated permutations $\pi^{:}(\mathcal I)$ of size $n$.
Under this bijection, the weak excedances of
$\pi^{:}(\mathcal I)$ are in positions $I_1$.
\end{lemma}


One particularly nice class of positroid cells is the \emph{TP}
or \emph{totally positive} Schubert cells.
\begin{definition}
A \emph{totally positive} Schubert cell is
the intersection of a Schubert cell with $\Grkn$.
\end{definition}

TP Schubert cells are indexed by $\Le$-diagrams such that all boxes are filled with
a $+$.


\section{A Deodhar decomposition of the Grassmannian}
\label{sec:project}

In this section we review Deodhar's decomposition of the 
flag variety  $G/B$ \cite{Deodhar}, and the parameterizations of the components due to Marsh 
and Rietsch \cite{MR}.  
We then give a Deodhar decomposition of the Grassmannian by projecting the usual
Deodhar decomposition of the flag variety to the Grassmannian.

\subsection{The flag variety}
In this paper we fix the group
$G = \SL_n = \SL_n(\R)$,
a maximal torus $T$, and opposite Borel subgroups $B^+$ and $B^-$, which
consist of the diagonal, upper-triangular, and lower-triangular
matrices, respectively.  
We let $U^+$ and $U^-$ be the unipotent radicals of $B^+$ and $B^-$; these
are the subgroups of upper-triangular and lower-triangular matrices with $1$'s on the diagonals.
For each $1 \leq i \leq n-1$ we have a homomorphism 
$\phi_i:{\rm SL}_2\to {\rm SL}_n$ such that
\[
\phi_i\begin{pmatrix} a& b\\c&d\end{pmatrix}=
\begin{pmatrix}
1 &             &       &       &           &     \\
   &\ddots  &        &      &            &        \\
   &             &   a   &   b  &          &       \\
   &             &   c    &   d  &         &       \\
   &            &          &       &  \ddots  &     \\
   &            &          &       &              & 1 
   \end{pmatrix} ~\in {\rm SL}_n,
\]
that is, $\phi_i$ replaces a $2\times 2$ block of the identity matrix with $\begin{pmatrix} a&b\\c&d\end{pmatrix}$.
Here $a$ is at the $(i+1)$st diagonal entry counting from the southeast corner.\footnote{Our
numbering differs from that in \cite{MR} in that
the rows of our matrices in $\SL_n$ are numbered from the bottom.}
We use this to construct $1$-parameter subgroups in $G$ 
(landing in $U^+$ and $U^-$,
respectively) defined by
\begin{equation*}
x_i(m) = \phi_i \left(
                   \begin{array}{cc}
                     1 & m \\ 0 & 1\\
                   \end{array} \right)  \text{ and }\ 
y_i(m) = \phi_i \left(
                   \begin{array}{cc}
                     1 & 0 \\ m & 1\\
                   \end{array} \right) ,\
\text{ where }m \in \R.
\end{equation*}

Let $W$ denote the Weyl group $ N_G(T) / T$, 
where $N_G(T)$ is the normalizer of $T$.  
The simple reflections $s_i \in W$ are given explicitly by
$s_i:= \dot{s_i} T$ where $\dot{s_i} :=
                 \phi_i \left(
                   \begin{array}{cc}
                     0 & -1 \\ 1 & 0\\
                   \end{array} \right)$
and any $w \in W$ can be expressed as a product $w = s_{i_1} s_{i_2}
\dots s_{i_m}$ with $m=\ell(w)$ factors.  We set $\dot{w} =
\dot{s}_{i_1} \dot{s}_{i_2} \dots \dot{s}_{i_m}$.
For $G = \SL_n$, we have $W = \Sym_n$, the symmetric group on $n$ letters,
and $s_i$ is the transposition exchanging $i$ and $i+1$.

We can identify the flag variety $G/B$ with the variety
$\B$ of Borel subgroups, via
\begin{equation*}
gB \longleftrightarrow g \cdot B^+ := gB^+ g^{-1}.
\end{equation*}
We have two opposite Bruhat
decompositions of $\B$:
\begin{equation*}
\mathcal B=\bigsqcup_{w\in W} B^+\dot w\cdot B^+=\bigsqcup_{v\in W}
B^-\dot v\cdot B^+.
\end{equation*}
We define the \emph{Richardson variety}
\begin{equation*}
\mathcal R_{v,w}:=B^+\dot w\cdot B^+\cap B^-\dot v\cdot B^+,
\end{equation*}
the intersection of opposite Bruhat cells. This intersection is empty
unless $v\le w$, in which case it is smooth of dimension
$\ell(w)-\ell(v)$, see \cite{KaLus:Hecke2,Lusztig2}.

\subsection{Distinguished expressions}

We now provide background on distinguished and positive
distinguished subexpressions, as in
\cite{Deodhar} and \cite{MR}. We will assume that the reader is familiar
with the (strong) Bruhat order $<$ on the Weyl group $W=\Sym_n$, and the 
basics of reduced expressions, as in  \cite{BB}.

Let $\w:= s_{i_1}\dots s_{i_m}$ be a reduced expression for $w\in W$.
A {\it subexpression} $\v$ of $\w$
is a word obtained from the reduced expression $\w$ by replacing some of
the factors with $1$. For example, consider a reduced expression in $\Sym_4$, say $s_3
s_2 s_1 s_3 s_2 s_3$.  Then $s_3 s_2\, 1\, s_3 s_2\, 1$ is a
subexpression of $s_3 s_2 s_1 s_3 s_2 s_3$.
Given a subexpression $\v$, 
we set $v_{(k)}$ to be the product of the leftmost $k$ 
factors of $\v$, if $k \geq 1$, and $v_{(0)}=1$.

\begin{definition}\label{d:Js}\cite{MR, Deodhar}
Given a subexpression $\v$ of a reduced expression $\w=
s_{i_1} s_{i_2} \dots s_{i_m}$, we define
\begin{align*}
J^{\circ}_\v &:=\{k\in\{1,\dotsc,m\}\ |\  v_{(k-1)}<v_{(k)}\},\\
J^{\Box}_\v\, &:=\{k\in\{1,\dotsc,m\}\ |\  v_{(k-1)}=v_{(k)}\},\\
J^{\bullet}_\v &:=\{k\in\{1,\dotsc,m\}\ |\  v_{(k-1)}>v_{(k)}\}.
\end{align*}
The expression  $\v$
is called {\it
non-decreasing} if $v_{(j-1)}\le v_{(j)}$ for all $j=1,\dotsc, m$,
e.g.\ $J^{\bullet}_\v=\emptyset$.
\end{definition}

\begin{definition}[Distinguished subexpressions]\cite[Definition 2.3]{Deodhar}
A subexpression $\v$ of $\w$ is called {\it distinguished}
if we have
\begin{equation}\label{e:dist}
v_{(j)}\le v_{(j-1)}\ s_{i_j}\qquad \text{for all
$~j\in\{1,\dotsc,m\}$}.
\end{equation}
In other words, if right multiplication by $s_{i_j}$ decreases the
length of $v_{(j-1)}$, then in a distinguished subexpression we
must have
$v_{(j)}=v_{(j-1)}s_{i_j}$.

We write $\v\prec\w$ if $\v$ is a distinguished subexpression of
$\w$.
\end{definition}

\begin{definition}[Positive distinguished subexpressions]
We call a
subexpression $\v$ of $\w$ a {\it positive distinguished subexpression} 
(or a PDS for short) if
 \begin{equation}\label{e:PositiveSubexpression}
v_{(j-1)}< v_{(j-1)}s_{i_j} \qquad \text{for all
$~j\in\{1,\dotsc,m\}$}.
 \end{equation}
In other words, it is distinguished and non-decreasing.
\end{definition}

\begin{lemma}\label{l:positive}\cite{MR}
Given $v\le w$ 
and a reduced expression $\w$ for $w$,
there is a unique PDS $\v_+$ for $v$ in $\w$.
\end{lemma}

\subsection{Deodhar components in the flag variety}

We now describe the Deodhar decomposition of the  flag variety.  This is a further
refinement of the decomposition of $G/B$ into Richardson varieties $\mathcal R_{v,w}$.
Marsh and Rietsch \cite{MR} gave explicit parameterizations for each Deodhar
component, identifying each one with a subset in the group.


\begin{definition}\cite[Definition 5.1]{MR}\label{d:factorization}
Let $\w=s_{i_1} \dots s_{i_m}$ be a reduced expression for $w$,
and let $\v$ be a distinguished subexpression.
Define a subset $G_{\v,\w}$ in $G$ by
\begin{equation}\label{e:Gvw}
G_{\v,\w}:=\left\{g= g_1 g_2\cdots g_m \left
|\begin{array}{ll}
 g_\ell= x_{i_\ell}(m_\ell)\dot s_{i_\ell}\inv& \text{ if $\ell\in J^{\bullet}_\v$,}\\
 g_\ell= y_{i_\ell}(p_\ell)& \text{ if $\ell\in J^{\Box}_\v$,}\\
 g_\ell=\dot s_{i_\ell}& \text{ if $\ell\in J^{\circ}_\v$,}
 \end{array}\quad \text{
for $p_\ell\in\R^*,\, m_\ell\in\R$. }\right. \right\}.
\end{equation}
There is an obvious map $(\R^*)^{|J^{\Box}_\v|}\times \R^{|J^{\bullet}_\v|}\to
G_{\v,\w}$ defined by the parameters $p_\ell$ and $m_\ell$ in
\eqref{e:Gvw}. 
For $v =w=1$ 
we define $G_{\v,\w}=\{1\}$.
\end{definition}

\begin{example}\label{ex:g}
Let $W=\Sym_5$, $\w = s_2 s_3 s_4 s_1 s_2 s_3$ and 
$\v = s_2 1 1 1 s_2 1$. 
Then the corresponding element $g\in G_{\v,\w}$ is given by $g=s_2 y_3(p_2)y_4(p_3)y_1(p_4)x_2(m_5)s_2^{-1}y_3(p_6)$,
which is 
\[
g=\begin{pmatrix}
1 & 0 & 0 & 0 & 0 \\
p_3 & 1 & 0 & 0 & 0 \\
0 & p_6 & 1 & 0 & 0 \\
p_2 p_3 & p_2-m_5 p_6 & -m_5 & 1 & 0 \\
0 & -p_4 p_6 & -p_4 & 0 & 1
\end{pmatrix}.
\]
\end{example}

The following result from \cite{MR} gives an explicit parametrization for
the Deodhar component $\mathcal R_{\v,\w}$.
We will take the description below as the \emph{definition}
of $\mathcal R_{\v,\w}$.

\begin{proposition}\label{p:parameterization}\cite[Proposition 5.2]{MR}
The map $(\R^*)^{|J^{\Box}_\v|}\times \R^{|J^{\bullet}_\v|}\to G_{\v,\w}$ from
Definition~\ref{d:factorization} is an isomorphism. The set
$G_{\v,\w}$ lies in $U^-\dot v\cap B^+\dot w B^+$, and the
assignment $g\mapsto g\cdot B^+$ defines an isomorphism
\begin{align}\label{e:parameterization}
G_{\v,\w}&~\overset\sim\To ~\mathcal R_{\v,\w}
\end{align}
between the subset $G_{\v,\w}$ of the group,
and the Deodhar component $\mathcal R_{\v,\w}$
 in the flag variety.
\end{proposition}

Suppose that for each $w\in W$ we choose a reduced expression 
$\w$ for $w$.  Then it follows from Deodhar's work (see \cite{Deodhar} and
\cite[Section 4.4]{MR}) that 
\begin{equation}\label{e:DeoDecomp}
\mathcal R_{v,w} = \bigsqcup_{\v \prec \w} \mathcal R_{\v,\w}\qquad  \text{ and }
\qquad  G/B=\bigsqcup_{w\in W}\left(\bigsqcup_{\v\prec \w} \mathcal
 R_{\v,\w}\right).
\end{equation}
These are called the \emph{Deodhar decompositions} of $\mathcal R_{v,w}$ and  $G/B$.

\begin{remark}\label{rem:KLpolys}
One may define the Richardson variety $\mathcal R_{v,w}$ over a finite field 
$\F_q$.  In this setting the number of points
determine the $R$-polynomials $R_{v,w}(q) = \#(\mathcal R_{v,w}(\F_q))$
introduced by Kazhdan and Lusztig \cite{KazLus} to give a formula
for the Kazhdan-Lusztig polynomials.  This was the original motivation
for Deodhar's work.
Therefore 
the isomorphisms
$\mathcal R_{\v,\w} \cong (\F_q^*)^{|J^{\Box}_\v} \times \F_q^{|J^{\bullet}_\v|}$
together with the decomposition \eqref{e:DeoDecomp} 
give formulas for the $R$-polynomials.
\end{remark}



\subsection{Projections of Deodhar components to the Grassmannian}\label{sec:projections}

Now we consider the projection of the Deodhar decomposition
to the Grassmannian $Gr_{k,n}$.
Let $W_k$ be the parabolic subgroup 
$\langle s_1,s_2,\dots,\hat{s}_{n-k},\dots,s_{n-1} \rangle$
of $W = \Sym_n$, and 
let $W^k$ denote the set of minimal-length
coset representatives of $W/W_k$. 
Recall that a \emph{descent} of a permutation $\pi$
is a position $j$ such that $\pi(j)>\pi(j+1)$.
Then $W^k$ is the subset of permutations of $\Sym_n$
which have at most one descent; and that descent must be in position $n-k$.

Let $\pi_k: G/B \to Gr_{k,n}$ be the projection from the flag variety to the Grassmannian.
For each $w \in W^k$ and $v \leq w$, define 
$\mathcal P_{v,w} = \pi_k(\mathcal R_{v,w})$.  Then 
by work of Lusztig \cite{Lusztig2},
$\pi_k$ is an isomorphism on $\mathcal P_{v,w}$, and 
we have a decomposition 
\begin{equation}\label{projected-Richardson}
Gr_{k,n} = \bigsqcup_{w \in W^k} \left(\bigsqcup_{v \leq w} \mathcal P_{v,w} \right).
\end{equation}

\begin{definition}
For each reduced decomposition $\w$ for $w \in W^k$, and each $\v \prec \w$,
we define the \emph{(projected) Deodhar component}
$\mathcal P_{\v,\w} = \pi_k(\mathcal R_{\v,\w})\subset Gr_{k,n}$.
\end{definition}

\begin{lemma}\label{rem:coincide}\cite[Remark 3.12, Lemma 3.13]{KW3}
The decomposition in \eqref{projected-Richardson} coincides with the
positroid stratification from Section \ref{sec:positroid}.  
The appropriate bijection 
between the strata is defined as follows.
Let $\Q^k$ denote
the set of pairs $(v,w)$ where $v\in W$, $w\in W^k$,  and $v\leq w$;
let $\Dec_n^k$ denote the set of decorated permutations in $S_n$
with $k$ weak excedances.
We consider both sets as partially ordered sets, where the cover
relation corresponds to containment of closures of the corresponding 
strata.  Then there is an order-preserving bijection 
$\Phi$ from $\Q^k$ to $\Dec_n^k$ which is defined as follows.
Let $(v,w) \in \Q^J$.  Then $\Phi(v,w) = (\pi, col)$
where $\pi = v w^{-1}$.  We also let $\pi^:(v,w)$ denote $\Phi(v,w)$.
To define $col$, we color any fixed point that occurs
in one of the positions $w(1), w(2), \dots, w(n-k)$
with the color $-1$, and color any other fixed point with the color $1$.
\end{lemma}

By Lemma \ref{rem:coincide},
$\mathcal P_{\v,\w}$ lies in the 
positroid stratum $S_{\pi^:}$.
Note that 
the strata $\mathcal P_{\v,\w}$ do not depend on the chosen
reduced decomposition
of $\w$  
\cite[Proposition 4.16]{KW2}.

Now if for each $w \in W^k$ we choose a reduced decomposition $\w$, then we have 
\begin{equation}\label{e:ProjDeoDecomp}
\mathcal P_{v,w} = \bigsqcup_{\v \prec \w} \mathcal P_{\v,\w}\qquad \text{ and }\qquad
 Gr_{k,n} =\bigsqcup_{w\in W^k} \left(\bigsqcup_{\v\prec \w} \mathcal
 P_{\v,\w}\right).
\end{equation}

Proposition \ref{p:parameterization} gives us a concrete way to 
think about the projected Deodhar components $\mathcal P_{\v,\w}$.  
The projection $\pi_k: G/B \to Gr_{k,n}$ maps 
each $g \in G_{\v,\w}$ to the span of its leftmost $k$ columns:
\[
g=\begin{pmatrix}
g_{n,n} & \dots& g_{n,n-k+1} & \dots & g_{n,1}  \\
\vdots & &\vdots&  & \vdots  \\
g_{1,n} & \dots& g_{1,n-k+1} & \dots & g_{1,1} \\
\end{pmatrix}
\quad \mapsto \quad
A=
\begin{pmatrix}
g_{1,n-k+1}&  \dots & g_{n,n-k+1}\\
\vdots &  & \vdots\\
g_{1,n} &  \dots & g_{n,n} \\
\end{pmatrix}.
\]
Alternatively, we may identify $A\in Gr_{k,n}$ with its image in the Pl\"ucker 
embedding. Let $e_i$ denote the column vector in $\R^n$ 
such that the $i$th entry from the bottom contains a $1$, and all other 
entries are $0$, e.g. $e_n=(1,0, \ldots,0)^T$, the transpose of the row vector $(1,0,\ldots,0)$.
Then the projection $\pi_k$ maps
each $g \in G_{\v,\w}$ (identified with $g \cdot B^+ \in \mathcal R_{\v,\w}$) to 
\begin{align}\label{Plucker}
g \cdot e_n \wedge \ldots \wedge e_{n-k+1} &=
\sum_{1\le 
j_1<\ldots<j_k\le n}\Delta_{j_1,\ldots,j_k}(A)e_{j_k}\wedge \cdots\wedge e_{j_1}.
\end{align}
That is, the Pl\"ucker coordinate $\Delta_{j_1,\ldots,j_k}(A)$ is given by
\[
\Delta_{j_1,\ldots,j_k}(A)=\langle e_{j_k}\wedge\cdots\wedge e_{j_1},\,
g\cdot e_n\wedge\cdots\wedge e_{n-k+1}\rangle,
\]
where $\langle\cdot,\cdot\rangle$ is the usual inner product on $\wedge^k\mathbb{R}^n$.

\begin{example}
We continue Example \ref{ex:g}.  Note that 
$w \in W^k$ where $k=2$.  Then the map $\pi_2: G_{\v,\w}\to Gr_{2,5}$ 
is given by
\[
g=\begin{pmatrix}
1 & 0 & 0 & 0 & 0 \\
p_3 & 1 & 0 & 0 & 0 \\
0 & p_6 & 1 & 0 & 0 \\
p_2 p_3 & p_2-m_5 p_6 & -m_5 & 1 & 0 \\
0 & -p_4 p_6 & -p_4 & 0 & 1
\end{pmatrix} \quad\longrightarrow \quad
A = 
\begin{pmatrix}
-p_4p_6&p_2-m_5 p_6 & p_6 & 1 & 0\\
0 & p_2 p_3 & 0 & p_3 & 1 \\
\end{pmatrix}.
\]
\end{example}



\section{Deodhar components in the Grassmannian and Go-diagrams}\label{Deodhar-combinatorics}
In this section we explain how to index the Deodhar components in
the Grassmannian $Gr_{k,n}$ by certain
tableaux called \emph{Go-diagrams}, which are 
fillings of 
Young diagrams by 
empty boxes, \raisebox{0.12cm}{\hskip0.14cm\circle*{7}\hskip-0.1cm}'s and \raisebox{0.12cm}{\hskip0.14cm\circle{7}\hskip-0.1cm}'s.
We refer to the symbols
 \raisebox{0.12cm}{\hskip0.14cm\circle*{7}\hskip-0.1cm} and \raisebox{0.12cm}{\hskip0.14cm\circle{7}\hskip-0.1cm}  as \emph{black} and \emph{white stones}.
Recall that $W_k = \langle s_1,s_2,\dots,\hat{s}_{n-k},\dots,s_{n-1} \rangle$
is a parabolic subgroup of $W =\Sym_n$ and $W^k$ is the set of minimal-length
coset representatives of $W/W_k$. 



We fix $k$ and $n$, and let $Q^k$ be the poset 
whose elements are the boxes in a 
$k \times (n-k)$ rectangle.  
If $b_1$ and $b_2$
are two adjacent boxes such that $b_2$ is immediately to the left or
immediately above $b_1$, we have a cover relation $b_1 \lessdot b_2$
in $Q^k$. The partial order on $Q^k$ is the transitive closure of
$\lessdot$. 
The middle and right diagram in the figure below show 
two linear extensions (or \emph{reading orders}) of 
the poset $Q^3$, where $n=8$.
Next we assign a labeling of the boxes of $Q^k$ by 
simple generators $s_i$, see the left diagram
of the figure below.  If $b$ is a box
labelled by $s_i$, we denote the simple generator labeling
$b$ by $s_b:= s_i$.

The following result can be found in \cite{Ste}.
\begin{proposition}
Fix $k$ and $n$.  The upper order ideals of $Q^k$
are in bijection with elements of $W^k$, and the different
reading orders of $Q^k$ allow us to compute all 
reduced expressions of elements of $W^k$.  More specifically,
let $Y$ be an upper order ideal of $Q^k$,
and choose a reading order $e$ for the
boxes of $Y$.  Then if we read the labels of $Y$
in the order specified by $e$, we will get a reduced
word for some $w \in W^k$, and this element $w$
does not depend on the choice of $e$. Therefore we may 
denote $Y$ by $O_w$.   Moreover,
if we let $e$ vary over all reading orders
for boxes of $Y$, we will
obtain all reduced expressions for $w$.
\end{proposition}

\setlength{\unitlength}{0.7mm}
\begin{center}
  \begin{picture}(60,30)
  
  \put(5,32){\line(1,0){45}}
  \put(5,23){\line(1,0){45}}
  \put(5,14){\line(1,0){45}}
  \put(5,5){\line(1,0){45}}
  \put(5,5){\line(0,1){27}}
  \put(14,5){\line(0,1){27}}
  \put(23,5){\line(0,1){27}}
  \put(32,5){\line(0,1){27}}
  \put(41,5){\line(0,1){27}}
  \put(50,5){\line(0,1){27}}

  \put(7,27){$s_5$}
  \put(16,27){$s_4$}
  \put(25,27){$s_3$}
  \put(34,27){$s_2$}
  \put(43,27){$s_1$}
  \put(7,18){$s_6$}
  \put(16,18){$s_5$}
  \put(25,18){$s_4$}
  \put(34,18){$s_3$}
  \put(43,18){$s_2$}
  \put(7,9){$s_7$}
  \put(16,9){$s_6$}
  \put(25,9){$s_5$}
  \put(34,9){$s_4$}
  \put(43,9){$s_3$}
 
  \end{picture} 
  \quad
  \begin{picture}(60,40)
  
  \put(5,32){\line(1,0){45}}
  \put(5,23){\line(1,0){45}}
  \put(5,14){\line(1,0){45}}
  \put(5,5){\line(1,0){45}}
  \put(5,5){\line(0,1){27}}
  \put(14,5){\line(0,1){27}}
  \put(23,5){\line(0,1){27}}
  \put(32,5){\line(0,1){27}}
  \put(41,5){\line(0,1){27}}
  \put(50,5){\line(0,1){27}}

  \put(7,26){$15$}
  \put(16,26){$14$}
  \put(25,26){$13$}
  \put(34,26){$12$}
  \put(43,26){$11$}
  \put(7,17){$10$}
  \put(16,17){$~9$}
  \put(25,17){$~8$}
  \put(34,17){$~7$}
  \put(43,17){$~6$}
  \put(7,8){$~5$}
  \put(16,8){$~4$}
  \put(25,8){$~3$}
  \put(34,8){$~2$}
  \put(43,8){$~1$}
  
 \end{picture} 
  \quad
  \begin{picture}(60,40)
  
  \put(5,32){\line(1,0){45}}
  \put(5,23){\line(1,0){45}}
  \put(5,14){\line(1,0){45}}
  \put(5,5){\line(1,0){45}}
  \put(5,5){\line(0,1){27}}
  \put(14,5){\line(0,1){27}}
  \put(23,5){\line(0,1){27}}
  \put(32,5){\line(0,1){27}}
  \put(41,5){\line(0,1){27}}
  \put(50,5){\line(0,1){27}}
  
  \put(7,26){$15$}
  \put(16,26){$12$}
  \put(25,26){$~9$}
  \put(34,26){$~6$}
  \put(43,26){$~3$}
  \put(7,17){$14$}
  \put(16,17){$11$}
  \put(25,17){$~8$}
  \put(34,17){$~5$}
  \put(43,17){$~2$}
  \put(7,8){$13$}
  \put(16,8){$10$}
  \put(25,8){$~7$}
  \put(34,8){$~4$}
  \put(43,8){$~1$}

  \end{picture} 
\end{center}

\begin{remark}
The upper order ideals of $Q^k$ can be identified with the Young diagrams
contained in a $k \times (n-k)$ rectangle (justified at the upper left),
and the reading orders
of $O_w$ can be identified with the \emph{reverse} standard tableaux of shape $O_w$, i.e. entries decrease from left to right in rows and from top to bottom
in columns.
\end{remark}

\subsection{Go-diagrams 
and labeled Go-diagrams}

Next we will identify distinguished subexpressions of reduced words 
for elements of $W^k$ with certain tableaux called \emph{Go-diagrams}.

\begin{definition}\label{def:Go}
Let $\w$ be a reduced expression for some $w \in W^k$,
and consider a distinguished subexpression $\v$ of $\w$.
Let $O_w$ be the Young diagram associated to $w$
and choose the reading order of its boxes which corresponds
to $\w$.  Then 
for each $k \in J_{\v}^{\circ}$ we will place a \raisebox{0.12cm}{\hskip0.15cm\circle{4}\hskip-0.15cm} in the corresponding box of $O_w$;
for each 
$k \in J_{\v}^{\bullet}$ we will place a \raisebox{0.12cm}{\hskip0.15cm\circle*{4}\hskip-0.15cm} in the corresponding box; 
and for each $k \in J_{\v}^{\Box}$ we will leave the corresponding box blank.
We call the resulting diagram a \emph{Go-diagram},
and refer to the symbols
\wstn ~ and \bstn~  
as {\it white} and {\it  black stones}.
\end{definition}

\begin{remark}
Note that a Go-diagram has no black stones if and only if 
$\v$ is a positive distinguished subexpression of $\w$.
In this case, if we replace the empty boxes by $+$'s 
and the white stones by $0$'s, we will get a 
$\Le$-diagram.  
Therefore,
slightly abusing terminology,
we will often refer to a Go-diagram with no black stones as a 
$\Le$-diagram.\footnote{Since
$\Le$-diagrams are a special case of Go-diagrams, one might also refer to them
as \emph{Lego} diagrams.}
Moreover, the Go-diagrams with no black
stones are in bijection with $\Le$-diagrams.  See 
\cite[Section 4]{KW3} for more details.
\end{remark}

\begin{example}\label{ex4-1}
Consider the upper order ideal $O_w$ which is $Q^k$ itself for $\Sym_5$ and $k=2$.
Then $Q^k$ is the poset shown in the left diagram.  Let us choose the reading order  indicated by the
labeling shown in the right diagram.
\setlength{\unitlength}{0.7mm}
\begin{center}
    \begin{picture}(50,25)
  
  \put(5,23){\line(1,0){27}}
  \put(5,14){\line(1,0){27}}
  \put(5,5){\line(1,0){27}}
  \put(5,5){\line(0,1){18}}
  \put(14,5){\line(0,1){18}}
  \put(23,5){\line(0,1){18}}
  \put(32,5){\line(0,1){18}}

  \put(7,18){$s_3$}
  \put(16,18){$s_2$}
  \put(25,18){$s_1$}
  \put(7,8){$s_4$}
  \put(16,8){$s_3$}
  \put(25,8){$s_2$}
 
  \end{picture} 
\qquad
  \begin{picture}(50,25)
  
  \put(5,23){\line(1,0){27}}
  \put(5,14){\line(1,0){27}}
  \put(5,5){\line(1,0){27}}
  \put(5,5){\line(0,1){18}}
  \put(14,5){\line(0,1){18}}
  \put(23,5){\line(0,1){18}}
  \put(32,5){\line(0,1){18}}

  \put(8,17){$6$}
  \put(17,17){$5$}
  \put(26,17){$4$}
  \put(8,8){$3$}
  \put(17,8){$2$}
  \put(26,8){$1$}
  
  \end{picture} 

\end{center}
\setlength{\unitlength}{0.7mm}

Now consider the distinguished
subexpressions $s_2 s_3 s_4 s_1 s_2 s_3$, $1 s_3 s_4 s_1 1 s_3$, and
$1 s_3 1 s_1 1 s_3$ of the 
reduced expression $\w = s_2 s_3 s_4 s_1 s_2 s_3$ corresponding
to our chosen reading order.
Among these three subexpressions, the 
first and second are
PDS's.  The corresponding Go-diagrams are as follows.
\begin{center}
    \begin{picture}(45,30)
  
  \put(5,25){\line(1,0){30}}
  \put(5,15){\line(1,0){30}}
  \put(5,5){\line(1,0){30}}
  \put(5,5){\line(0,1){20}}
  \put(15,5){\line(0,1){20}}
  \put(25,5){\line(0,1){20}}
  \put(35,5){\line(0,1){20}}

  \put(8,20){\hskip0.15cm\circle{5}}
  \put(18,20){\hskip0.15cm\circle{5}}
  \put(28,20){\hskip0.15cm\circle{5}}
  \put(8,10){\hskip0.15cm\circle{5}}
  \put(18,10){\hskip0.15cm\circle{5}}
  \put(28,10){\hskip0.15cm\circle{5}}
 
  \end{picture} 
\qquad
     \begin{picture}(45,30)
  
  \put(5,25){\line(1,0){30}}
  \put(5,15){\line(1,0){30}}
  \put(5,5){\line(1,0){30}}
  \put(5,5){\line(0,1){20}}
  \put(15,5){\line(0,1){20}}
  \put(25,5){\line(0,1){20}}
  \put(35,5){\line(0,1){20}}

  \put(8,20){\hskip0.15cm\circle{5}}
  \put(18,20){}
  \put(28,20){\hskip0.15cm\circle{5}}
  \put(8,10){\hskip0.15cm\circle{5}}
  \put(18,10){\hskip0.15cm\circle{5}}
  \put(28,10){}
 
  \end{picture} 
\quad
    \begin{picture}(45,30)
  
  \put(5,25){\line(1,0){30}}
  \put(5,15){\line(1,0){30}}
  \put(5,5){\line(1,0){30}}
  \put(5,5){\line(0,1){20}}
  \put(15,5){\line(0,1){20}}
  \put(25,5){\line(0,1){20}}
  \put(35,5){\line(0,1){20}}

  \put(8,20){\hskip0.15cm\circle*{5}}
  \put(18,20){}
  \put(28,20){\hskip0.15cm\circle{5}}
  \put(8,10){}
  \put(18,10){\hskip0.15cm\circle{5}}
  \put(28,10){}
 
  \end{picture} 
\end{center}
\end{example}

The reader might worry that Definition \ref{def:Go} 
has too much dependence on the choice of reduced expression
$\w$, or equivalently on the reading order $e$.
However, we have the following result.

\begin{proposition}\cite[Proposition 4.5]{KW3}
Let $D$ be a Go-diagram in a Young diagram $O_w$.
Choose any reading order $e$ for the boxes of $O_w$.
Let $\w$ be the corresponding reduced expression
for $w$.  Let $\v(D)$ be the subexpression of $\w$
obtained as follows: if a box $b$ of $D$ contains a black
or white stone then the corresponding simple generator 
$s_b$ is present in the subexpression, while if $b$
is empty, we omit the corresonding simple generator.
Then we have the following:
\begin{enumerate}
\item
the element $v:=v(D)$ is independent of the choice of reading
word $e$.
\item
whether $\v(D)$ is a PDS depends only on $D$ (and not $e$).
\item whether $\v(D)$ is distinguished depends only on $D$ (and not on $e$).
\end{enumerate}
\end{proposition}

\begin{definition}\label{def:DtoPi}
Let $O_w$ be an upper order ideal of $Q^k$, where 
$w\in W^k$ and $W = S_n$.
Consider a Go-diagram $D$ of shape $O_w$;
this is contained in a $k \times (n-k)$ rectangle,
and the shape $O_w$ gives rise to a lattice path
from the northeast corner to the 
southwest corner of the rectangle.  Label the steps
of that lattice path from $1$ to $n$; this gives
a natural labeling to every row and column of the rectangle.
We now let $v:=v(D)$, and
we define $\pi^:(D)$ to be the decorated permutation
$(\pi(D),col)$ where $\pi = \pi(D) = vw^{-1}$.  The fixed points of $\pi$
correspond precisely to rows and columns of the rectangle with no $+$'s.
If there are no $+$'s in the row (respectively, column) labeled by $h$,
then $\pi(h)=h$ and this fixed point gets colored with color $1$ (respectively,
$-1$.)
\end{definition}

\begin{remark}
It follows from Lemma \ref{rem:coincide} 
that the projected Deodhar component $\mathcal P_D$ corresponding to 
$D$ is contained in the positroid stratum $S_{\pi^:(D)}$.
\end{remark}

\begin{remark}\label{rem:KLpolys2}
Recall from Remark \ref{rem:KLpolys}
that the isomorphisms
$\mathcal R_{\v,\w} \cong (\F_q^*)^{|J^{\Box}_\v} \times \F_q^{|J^{\bullet}_\v|}$
together with the decomposition \eqref{e:DeoDecomp} 
give formulas for the $R$-polynomials.  Therefore a good 
characterization of Go-diagrams
could lead to explicit formulas for the corresponding $R$-polynomials.
\end{remark}

If we choose a reading order of $O_w$, then we will also associate
to a Go-diagram of shape $O_w$ a \emph{labeled Go-diagram},
as defined below.  Equivalently,
a labeled Go-diagram is associated to a pair $(\v,\w)$.

\begin{definition}\label{def:pi}
Given a reading order of $O_w$ and a Go-diagram of shape $O_w$,
we obtain a \emph{labeled Go-diagram} by replacing each \raisebox{0.12cm}{\hskip0.14cm\circle{4}\hskip-0.15cm} with a $1$,
each \raisebox{0.12cm}{\hskip0.14cm\circle*{4}\hskip-0.15cm} with a $-1$, and putting a $p_i$ in each blank square $b$,
where the subscript $i$ corresponds to the label of $b$ inherited from the 
reading order.
\end{definition}

The labeled Go-diagrams corresponding to 
the Go-diagrams  from
Example \ref{ex4-1} are:
\setlength{\unitlength}{0.7mm}
\begin{center}
    \begin{picture}(45,30)
  
  \put(5,25){\line(1,0){30}}
  \put(5,15){\line(1,0){30}}
  \put(5,5){\line(1,0){30}}
  \put(5,5){\line(0,1){20}}
  \put(15,5){\line(0,1){20}}
  \put(25,5){\line(0,1){20}}
  \put(35,5){\line(0,1){20}}

  \put(8,19){1}
  \put(18,19){1}
  \put(28,19){1}
  \put(8,9){1}
  \put(18,9){1}
  \put(28,9){1}
 
  \end{picture} 
\qquad
     \begin{picture}(45,30)
  
  \put(5,25){\line(1,0){30}}
  \put(5,15){\line(1,0){30}}
  \put(5,5){\line(1,0){30}}
  \put(5,5){\line(0,1){20}}
  \put(15,5){\line(0,1){20}}
  \put(25,5){\line(0,1){20}}
  \put(35,5){\line(0,1){20}}

  \put(8,19){1}
  \put(18,19){$p_5$}
  \put(28,19){$1$}
  \put(8,9){$1$}
  \put(18,9){$1$}
  \put(28,9){$p_1$}
 
  \end{picture} 
\quad
    \begin{picture}(45,30)
  
  \put(5,25){\line(1,0){30}}
  \put(5,15){\line(1,0){30}}
  \put(5,5){\line(1,0){30}}
  \put(5,5){\line(0,1){20}}
  \put(15,5){\line(0,1){20}}
  \put(25,5){\line(0,1){20}}
  \put(35,5){\line(0,1){20}}

  \put(7,19){$-1$}
  \put(18,19){$p_5$}
  \put(28,19){$1$}
  \put(8,9){$p_3$}
  \put(18,9){$1$}
  \put(28,9){$p_1$}
 
  \end{picture} 
\end{center}

\subsection{Pl\"ucker coordinates for  projected Deodhar components}\label{sec:positivitytest}

Consider $\mathcal P_{\v,\w} \subset Gr_{k,n}$, where 
$\w$ is a reduced expression for $w \in W^k$ and $\v \prec \w$.
Here we provide some formulas for the Pl\"ucker coordinates
of the elements of $\mathcal P_{\v,\w}$, in terms of the parameters used
to define $G_{\v,\w}$.  Some of these formulas are related to corresponding
formulas for $G/B$ in \cite[Section 7]{MR}.  

\begin{theorem}\label{th:maxmin}\cite[Lemma 5.1, Theorem 5.2]{KW3}
Choose any element $A$ of $\mathcal P_{\v,\w} \subset Gr_{k,n}$, 
in other words, 
$A = \pi_k(g)$ for some $g \in G_{\v,\w}$.
Then the lexicographically minimal and maximal 
nonzero Pl\"ucker coordinates of $A$ are 
$\Delta_I$ and $\Delta_{I'}$, where 
\begin{equation*}
I = w\,\{n,n-1,\dots,n-k+1\} \qquad \text{ and }\qquad
I' = v\,\{n,n-1,\dots,n-k+1\}.
\end{equation*}

Moreover, if we write $g = g_1 \dots g_m$ as in 
Definition \ref{d:factorization}, then
\begin{equation}
\Delta_I(A) = (-1)^{|J_{\v}^{\bullet}|} \prod_{i\in J_{\v}^{\Box}} p_i\qquad  \text{ and }\qquad
\Delta_{I'}(A) = 1.
\end{equation}
Note that $\Delta_I(A)$ equals the product of all the labels
from the labeled Go-diagram associated to $(\v,\w)$.
\end{theorem}

This theorem can be extended to provide a formula for some other Pl\"ucker coordinates.
Let $b$ be any box of $D$.  We can choose a linear extension $e$
of the boxes of $D$ which orders the boxes which are weakly southeast of 
$b$ (the \emph{inner boxes}) before the rest (the \emph{outer boxes}).  This gives rise to a reduced expression 
$\w$ and subexpression $\v$, as well as reduced expressions 
$\w_b^{\In} = \w^{\In}$, $\v_b^{\In} = \v^{\In}$, 
$\w_b^{\Out} = \w^{\Out}$, and $\v_b^{\Out} = \v^{\Out}$, which are  
obtained by restricting $\w$ and $\v$
to the inner and outer boxes, respectively.

\begin{theorem}\label{p:Plucker}\cite[Theorem 5.6]{KW3}
Let $\w=s_{i_1} \dots s_{i_m}$ be a reduced expression for $w \in W^k$ and $\v \prec \w$,
and let $D$ be the corresponding Go-diagram.  Choose any box $b$ of $D$,
and let $v^{\In} = v_b^{\In}$ and $w^{\In} = w_b^{\In}$, and
 $v^{\Out} = v_b^{\Out}$ and $w^{\Out} = w_b^{\Out}$.
Let $A = \pi_k(g)$ for any $g \in G_{\v,\w}$, and let
$I = w\{n,n-1,\dots,n-k+1\}$.
Define $I_b = v^{\In} (w^{\In})^{-1} I \in {[n] \choose k}$.
If we write $g = g_1 \dots g_m$ as in Definition \ref{d:factorization}, then
\begin{equation}
\Delta_{I_b}(A) = (-1)^{|J_{\v^{\Out}}^{\bullet}|} \prod_{i\in J_{\v^{\Out}}^\Box} p_i. 
\end{equation}
Note that $\Delta_{I_b}(A)$ equals the product of all the labels
in the ``out" boxes of the labeled Go-diagram.
\end{theorem}

\begin{example}
Consider $A\in Gr_{3,7}$ which is the projection of $g\in G_{\v,\w}$ with
\[
\w=s_3s_4s_5s_6s_2s_3s_4s_5s_1s_2s_3s_4,\qquad \v=s_3\,1\,s_5\,1\,s_2s_3\,1\,s_5\,1\,s_2s_3s_4.
\]
Then $v=s_2s_4$ and $\pi=vw^{-1}=(4,6,7,1,3,2,5)$.
The Go-diagram and the labeled Go-diagram are as follows:
\medskip
\[
\young[4,4,4][10][\hskip0.4cm\circle{5},\hskip0.4cm\circle*{5},\hskip0.4cm\circle*{5},,\hskip0.4cm\circle*{5},,\hskip0.4cm\circle{5},\hskip0.4cm\circle{5},,\hskip0.4cm\circle{5},,
\hskip0.4cm\circle{5}]\hskip1.5cm \young[4,4,4][10][$1$,$-1$,$-1$,$p_9$,$-1$,
$p_{7}$,$1$,$1$,$p_4$,$1$,$p_2$,$1$]
\hskip1.5cm\young[4,4,4][10][467,456,245,234,167,156,125,123,127,125,125,123]
\medskip
\]
\medskip
\noindent
Each box $b$ in the diagram on the right 
show the subset $I_b$.
Then we have 
\begin{align*}
\Delta_{1,2,3}=-p_2p_4p_7p_9,~\Delta_{1,2,5}=-p_4p_7p_9,~\Delta_{1,2,7}=-p_7p_9,~\Delta_{1,5,6}=-p_4p_9,\\
\Delta_{1,6,7}=p_9,~\Delta_{2,3,4}=-p_2p_4p_7,~\Delta_{2,4,5}=p_4p_7,~\Delta_{4,5,6}=-p_4,~\Delta_{4,6,7}=1.
\end{align*}

\end{example}


\section{Soliton solutions to the KP equation
and their contour plots}\label{soliton-background}

We now explain how to obtain a soliton solution 
to the KP equation from a point of $Gr_{k,n}$.

\subsection{From a point of the Grassmannian to a $\tau$-function.}

We start by fixing real parameters $\kappa_i$ such that 
$\kappa_1~<~\kappa_2~\cdots~<\kappa_n$,
which  are {\it generic}, in the sense
that 
the sums $\sum_{j=1} ^p\kappa_{i_j}$ are all distinct for any $p$ with
$1<p<n$.  We also 
assume that the differences between consecutive $\kappa_i$'s are similar,
that is, 
$\kappa_{i+1}-\kappa_i$ is of order one.

Let $\{E_i;i=1,\ldots,m\}$ be a set of exponential functions in $(x,y,t)\in \mathbb{R}^3$ defined by
\[
E_i(x,y,t):=\exp \theta_i(x,y,t) \quad\text{ where }\quad
\theta_i(x,y,t) = \kappa_i x + \kappa_i^2 y + \kappa_i^3 t.
\]
Thinking of each $E_i$ as a function of $x$,
we see that the elements $\{E_i\}$ are linearly independent,
because their Wronskian determinant with
respect to $x$ is non zero: 
\[
{\rm Wr}(E_1,\ldots,E_n):=\det [(E_i^{(j-1)})_{1\le i,j\le n}]=\prod_{i<j}(\kappa_j-\kappa_i)\,E_1\cdots E_n\ne 0.
\]
Here $E_i^{(j)}:=\partial^jE_i/\partial x^j=\kappa_i^jE_i$.

Let $A$ be a full rank $k \times n$ matrix.  We define 
a set of functions $\{f_1,\ldots,f_k\}$ by 
\[
(f_1,f_2,\ldots,f_k)^T = A\cdot (E_1,E_2,\ldots, E_n)^T,
\]
where $(\ldots)^T$ denotes the transpose of the vector $(\ldots)$.

Since the  exponential functions
 $\{E_i\}$ are linearly independent, we identify them as a basis of $\mathbb{R}^n$, and then
$\{f_1,\ldots,f_k\}$ spans a $k$-dimensional subspace.
This identification can be seen, more precisely, as $E_i\leftrightarrow (1,\kappa_i,\ldots,\kappa_i^{n-1})^T\in \mathbb{R}^n$.
This subspace  depends only on which point of 
the Grassmannian $Gr_{k,n}$
the matrix $A$ represents, so we can identify 
the space of subspaces 
$\{f_1,\ldots,f_k\}$ with $Gr_{k,n}$.

The {\it $\tau$-function} of $A$ is defined by
\begin{equation}\label{tauA}
\tau_A(x,y,t):={\rm Wr}(f_1,f_2,\ldots,f_k).
\end{equation}
For $I=\{i_1,\dots,i_k\} \in {[n] \choose k}$, we set
\[
E_I(x,y,t):={\rm Wr}(E_{i_1},E_{i_2},\ldots, E_{i_k})=\prod_{\ell<m}(\kappa_{i_m}-\kappa_{i_{\ell}})\,E_{i_1}\cdots E_{i_k}\,>0.
\]
Applying the Binet-Cauchy identity to the fact that  $f_j=\sum_{i=1}^n a_{ji}E_i$, 
we get \begin{equation}\label{tau}
\tau_A(x,y,t)=\sum_{I\in\binom{[n]}{k}}\Delta_I(A)\,E_I(x,y,t).
\end{equation}
It follows that if  $A\in \Grkn$, then 
$\tau_A>0$ for all $(x,y,t)\in\mathbb{R}^3$.  

Thinking of $\tau_A$ as a function of $A$, we note 
from \eqref{tau} that
the $\tau$-function encodes the information of the 
Pl\"ucker
embedding.  More specifically, 
if we identify each function $E_I$
with $I=\{i_1,\ldots,i_k\}$  
with the wedge product $E_{i_1}\wedge\cdots\wedge E_{i_k}$ (recall the identification
$E_i\leftrightarrow (1,\kappa_i,\ldots,\kappa_i^{n-1})^T$), 
then the map  $\tau:Gr_{k,n}\hookrightarrow \mathbb{RP}^{\binom{n}{k}-1}$,
$A\mapsto \tau_A$ has the Pl\"ucker coordinates as coefficients.

\subsection{From the $\tau$-function to solutions of the KP equation}
The KP equation for $u(x,y,t)$ 
\[
\frac{\partial}{\partial x}\left(-4\frac{\partial u}{\partial t}+6u\frac{\partial u}{\partial x}+\frac{\partial^3u}{\partial x^3}\right)+3\frac{\partial^2 u}{\partial y^2}=0
\]
was proposed by Kadomtsev and Petviashvili in 1970 \cite{KP70}, in order to 
study the stability of the soliton solutions of the Korteweg-de Vries (KdV) equation
under the influence of weak transverse perturbations.
The KP equation can be also used to describe  two-dimensional shallow 
water wave phenomena (see for example \cite{K10}).   This equation is now considered to be a prototype of
an integrable nonlinear partial differential equation.
For more background, see \cite{NMPZ84, D91, AC91, H04, MJD00}.  

It is well known
(see \cite{H04, CK1, CK2, CK3})   
that 
the $\tau$-function defined in \eqref{tauA} 
provides a soliton solution of the KP equation,
\begin{equation}\label{KPsolution}
u_A(x,y,t)=2\frac{\partial^2}{\partial x^2}\ln\tau_A(x,y,t).
\end{equation}

It is easy to show that if $A \in \Grkn$, then 
such a solution $u_A(x,y,t)$ is regular for all $(x,y,t)\in\mathbb{R}^3$.
For this reason we are interested in solutions $u_A(x,y,t)$
of the KP equation which come from points $A$ of 
$\Grkn$.
Throughout this paper when we speak of a {\it soliton solution to the KP equation},
we will mean a solution $u_A(x,y,t)$ which has form
\eqref{KPsolution}.


\subsection{Contour plots of soliton solutions}
\label{sec:solgraph}

\begin{figure}[h]
\begin{center}
\includegraphics[height=1.7in]{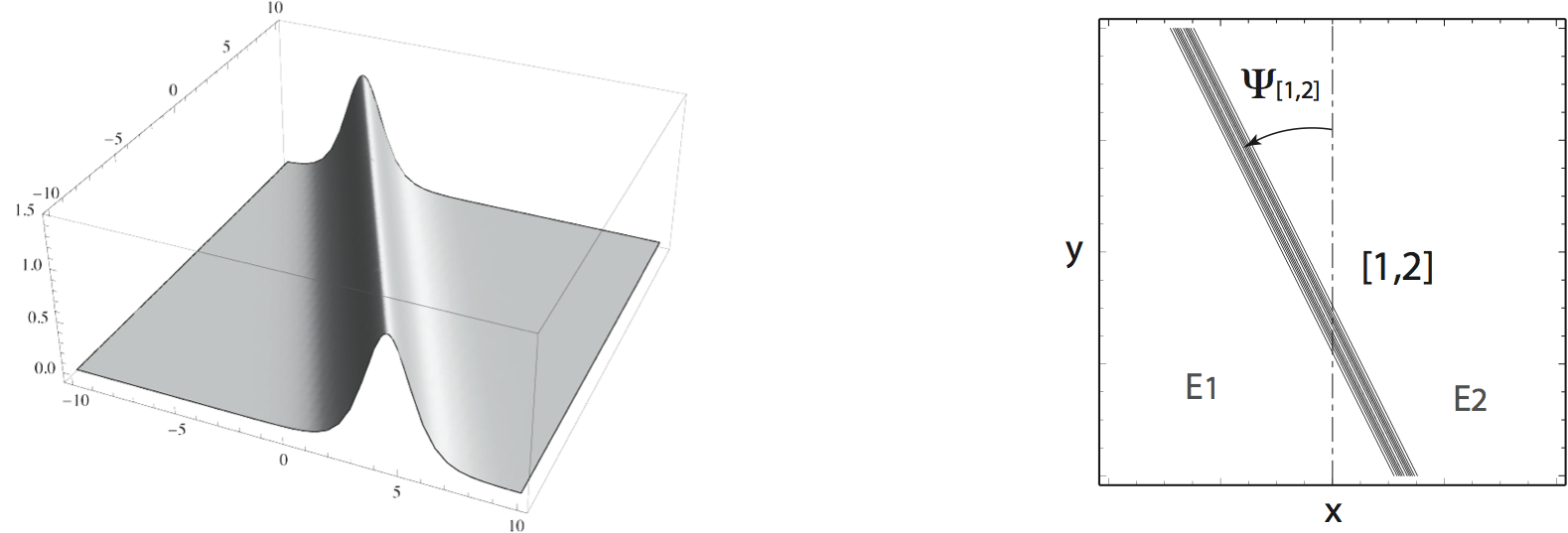}
\par
\end{center}
\caption{A line-soliton solution $u_A(x,y,t)$ where $A=(1,1) \in (Gr_{1,2})_{\geq 0}$, depicted via
the 3-dimensional profile
$u_A(x,y,t)$, and the level sets of $u_A(x,y,t)$ for some $t$.
$E_i$ represents the dominant exponential in each region.
\label{fig:1soliton}}
\end{figure}

One can visualize a solution $u_A(x,y,t)$ to the KP equation
by drawing 
level sets of the solution in the $xy$-plane, when the coordinate $t$ is fixed.
For each $r\in \mathbb{R}$, we denote the corresponding level set by 
\[
C_r(t):=\{(x,y)\in\mathbb{R}^2: u_A(x,y,t)=r\}. 
\]
Figure \ref{fig:1soliton} depicts both a three-dimensional image of a solution $u_A(x,y,t)$,
as well as multiple level sets $C_r$. 
These level sets are lines parallel to the line of the wave peak.

To study the behavior of $u_A(x,y,t)$ for $A\in S_{\M} \subset Gr_{k,n}$,
we consider the dominant exponentials in the $\tau$-function \eqref{tau} at each point $(x,y,t)$.  First we write the $\tau$-function in the form
\begin{align*}
\tau_A(x,y,t)&=\sum_{J\in\binom{[n]}{k}}\Delta_J(A)E_J(x,y,t)\\
&=\sum_{J\in\mathcal{M}}
\exp\left(\sum_{i=1}^n(\kappa_{j_i}x+\kappa_{j_i}^2y+\kappa_{j_i}^3t)+\ln(\Delta_J(A)K_J)\right), 
\end{align*}
where $K_J:=\prod_{\ell < m} (\kappa_{j_m}-\kappa_{j_{\ell}})>0$.
Note that in general the terms $\ln(\Delta_J(A)K_J)$ could be imaginary
when some $\Delta_J(A)$ are negative.

Since we are interested in the behavior of the soliton solutions when the variables
$(x,y,t)$ are on a large scale, we rescale the variables with a small positive number $\epsilon$, 
\[
x~\longrightarrow ~\frac{x}{\epsilon},\qquad y~\longrightarrow~\frac{y}{\epsilon},\qquad
t~\longrightarrow~\frac{t}{\epsilon}.
\]
This leads to
\[
\tau_A^{\epsilon}(x,y,t)
=\sum_{J\in\mathcal{M}}
\exp\left(\frac{1}{\epsilon}\,\sum_{i=1}^n(\kappa_{j_i}x+\kappa_{j_i}^2y+\kappa_{j_i}^3t)+\ln(\Delta_J(A)K_J)\right).
\]
Then we define a function $f_A(x,y,t)$ as the limit
\begin{equation}\label{f-contour}
\begin{array}{lll}
{f}_A(x,y,t)&=\displaystyle{\lim_{\epsilon\to 0}\epsilon\ln\left(\tau^{\epsilon}_A(x,y,t)\right)}\\[1.5ex]
&=  \underset{J\in\mathcal{M}}\max \left\{
                     \sum_{i=1}^k (\kappa_{j_i} x + \kappa^2_{j_i} y +\kappa^3_{j_i} t)\right\}.
\end{array}
\end{equation}
Since the above function depends only on the collection $\mathcal{M}$,
we also denote it as $f_{\mathcal{M}}(x,y,t)$.

\begin{definition}\label{contour}
Given a solution $u_A(x,y,t)$ of the KP equation as in 
(\ref{KPsolution}), we define its \emph{contour plot} $\mathcal{C}(u_A)$ to be the 
locus in $\mathbb{R}^3$ where $f_A(x,y,t)$ is not linear.  If we fix $t=t_0$, 
then we let $\mathcal{C}_{t_0}(u_A)$ be the locus in $\mathbb{R}^2$ where $f_A(x,y,t=t_0)$ is not 
linear, and we also refer to this as a \emph{contour plot}.
Because these contour plots depend only on $\mathcal{M}$ and not on $A$, we also refer to them
as $\mathcal{C}(\mathcal{M})$ and $\mathcal{C}_{t_0}(\mathcal{M})$.
\end{definition}

\begin{remark}
The contour plot approximates the locus where 
$|u_A(x,y,t)|$ takes on its maximum values or is singular.
\end{remark}

\begin{remark}
Note that the contour plot generated by the function $f_A(x,y,t)$ at $t=0$ consists of a set of semi-infinite lines attached to the origin $(0,0)$ in the $xy$-plane.
And if $t_1$ and $t_2$ have the same sign,
then the corresponding contour plots
$\CC_{t_1}(\M)$ and $\CC_{t_2}(\M)$ are self-similar.

Also note that because  our definition of the contour plot ignores the constant 
terms $\ln(\Delta_J(A)K_J)$, there are no phase-shifts in our picture,
and the contour plot for $f_A(x,y,t) = f_{\M}(x,y,t)$ does
not depend on the signs of the Pl\"ucker coordinates.
\end{remark}

It follows from Definition \ref{contour} that $\CC(u_A)$ 
and $\CC_{t_0}(u_A)$
are piecewise linear subsets of $\R^3$ and $\R^2$, respectively, of 
codimension $1$.  In fact it is easy to verify the following.

\begin{proposition}\cite[Proposition 4.3]{KW2}
If each $\kappa_{i}$ is an integer, then 
$\CC(u_A)$ is a tropical hypersurface in $\R^3$,
and $\CC_{t_0}(u_A)$ is a tropical hypersurface (i.e. a tropical curve)
in $\R^2$.  
\end{proposition}

The contour plot $\CC_{t_0}(u_A)$ consists of line segments 
called \emph{line-solitons}, some of 
which have finite length, while others are unbounded
and extend in the $y$ direction to 
$\pm \infty$. 
Each region of the complement of 
$\CC_{t_0}(u_A)$ 
in $\R^2$ is a domain of linearity for $f_A(x,y,t)$, 
and hence each region is  
naturally associated to a {\it dominant exponential}  
$\Delta_J(A) E_J(x,y,t)$ from the $\tau$-function \eqref{tau}.
We label this region by $J$ or $E_J$.
Each line-soliton
represents a balance between two dominant exponentials
in the 
$\tau$-function.



Because of the genericity of the $\kappa$-parameters, the following
lemma is immediate.
\begin{lemma}\label{separating}\cite[Proposition 5]{CK3}
The index sets of the 
dominant exponentials of the $\tau$-function in adjacent regions
of the contour plot in the $xy$-plane are of the form $\{i,l_2,\dots,l_k\}$ and
$\{j, l_2,\dots,l_k\}$.
\end{lemma}

We call the line-soliton separating the two dominant exponentials in 
Lemma \ref{separating} a {\it line-soliton of type $[i,j]$}.  Its
equation is 
\begin{equation}\label{eq-soliton}
x+(\kappa_i+\kappa_j)y+(\kappa_i^2+\kappa_i\kappa_j+\kappa_j^2)t=0.
\end{equation}

\begin{remark}\label{slope}
Consider a line-soliton given by (\ref{eq-soliton}).
Compute the angle $\Psi_{[i,j]}$ between 
the positive $y$-axis and 
the line-soliton 
of type $[i,j]$, 
measured in the counterclockwise direction, so that the negative $x$-axis
has an angle of $\frac{\pi}{2}$ and the positive $x$-axis has an 
angle of $-\frac{\pi}{2}$. Then $\tan \Psi_{[i,j]} = \kappa_i+\kappa_j$,
so we refer to $\kappa_i+\kappa_j$ as the \emph{slope} of the 
$[i,j]$ line-soliton (see Figure \ref{fig:1soliton}).
\end{remark}

In Section \ref{sec:solitonplabic} we will explore 
the combinatorial structure of contour plots,
that is, the ways in which line-solitons may interact. 
Generically we expect a point at which several line-solitons
meet to have degree $3$; we regard such a point as a trivalent 
vertex.  Three line-solitons meeting at a trivalent vertex
exhibit a {\it resonant interaction} (this corresponds to the 
{\it balancing condition} for a tropical curve). 
See \cite[Section 4.2]{KW2}.  
One may also have two line-solitons which cross over
each other, forming an $X$-shape: we call this
an \emph{$X$-crossing}, but do not regard it as a vertex.
See Figure \ref{contour-soliton}. 
Vertices of degree greater than $4$
are also possible.  

\begin{definition}\label{def:blackwhiteX}
Let $i<j<k<\ell$ be positive integers.  
An $X$-crossing involving two line-solitons of types
$[i,k]$ and $[j,\ell]$ is called a 
\emph{black $X$-crossing}.  An $X$-crossing involving
two line-solitons of types $[i,j]$ and $[k,\ell]$, or 
of types $[i,\ell]$ and $[j,k]$, is called a 
\emph{white $X$-crossing}.
\end{definition}

\begin{definition}\label{def:generic}
A contour plot $\CC_{t}(u_A)$ is called \emph{generic} if all interactions of line-solitons
are at trivalent vertices or are $X$-crossings. 
\end{definition}

\section{Unbounded line-solitons at $y\gg0$ and $y\ll0$}\label{sec:unbounded}

In this section we explain that the unbounded line-solitons
at $|y|\gg0$ of a contour plot
$\CC_{t}(u_A)$ are determined by which positroid stratum
contains  $A$. 
Conversely, the unbounded line-solitons of
$\CC_{t}(u_A)$ determine which positroid stratum
$A$ lies in.

\begin{theorem}\label{thm:soliton-perm}\cite[Theorem 8.1]{KW3}
Let $A \in Gr_{k,n}$ lie in the positroid stratum 
$S_{\pi^:}$, where $\pi^: = 
(\pi,col)$.
Consider the contour plot
$\CC_{t}(u_A)$ for any time $t$.  
Then the excedances (respectively, nonexcedances)
of $\pi$ are in bijection with the 
unbounded line-solitons of $\CC_t(u_A)$ at $y\gg0$ (respectively, $y\ll0$).
More specifically, in 
$\CC_t(u_A)$, 
\begin{itemize}
\item[(a)]  
there is 
an unbounded line-soliton of  $[i,h]$-type at $y\gg0$
if and only if
$\pi(i)=h$ for $i<h$, 
\item[(b)] 
there is 
an unbounded line-soliton of $[i,h]$-type at $y\ll 0$  
if and only if
$\pi(h)=i$ for $i<h$.
\end{itemize}
Therefore $\pi^:$ determines
the unbounded line-solitons at $y\gg 0$ and $y\ll0$
of 
$\CC_{t}(u_A)$ for any time $t$.  

Conversely, given a 
contour plot
$\CC_{t}(u_A)$ at any time $t$ where $A \in Gr_{k,n}$,
one can construct $\pi^:=(\pi,col)$ such that $A \in S_{\pi^:}$ as follows.
The excedances and nonexcedances of $\pi$ are constructed as above
from the unbounded line-solitons. 
If there is an $h \in [n]$ such that $h\in J$ for every 
dominant exponential $E_J$ labeling the contour plot, 
then set $\pi(h)=h$ with $col(h)=1$.
If there is an $h\in [n]$ such that $h\notin J$ for any
dominant exponential $E_J$ labeling the contour plot,
then set $\pi(h)=h$ with $col(h)=-1$.
\end{theorem}

\begin{remark} Chakravarty and Kodama \cite[Prop. 2.6 and 2.9]{CK1} and 
\cite[Theorem 5]{CK3} associated a derangement to each 
\emph{irreducible} element $A$
in the \emph{totally non-negative part} $(Gr_{k,n})_{\geq 0}$ 
of the Grassmannian. 
Theorem \ref{thm:soliton-perm} generalizes their result by 
dropping the hypothesis of irreducibility and extending the 
setting from 
$(Gr_{k,n})_{\geq 0}$ to $Gr_{k,n}$.
\end{remark}

\begin{example}\label{ex:Gr49CP}
Consider some $A\in Gr_{4,9}$ which is the projection of an element
 $g\in G_{\v,\w}$ with
\begin{equation*}
{\bf w}=s_7s_8s_4s_5s_6s_7s_2s_4s_5s_6s_1s_2s_3s_4s_5 \quad \text{ and }\quad
{\bf v}=s_711s_51s_7s_21s_4111s_21s_4s_5.
\end{equation*}
Then  $v=1$ and $\pi=vw^{-1} = 
(6,7,1,8,2,3,9,4,5).$
The matrix $g\in G_{\v,\w}$ is given by
\begin{align*}
g=
\dot s_7y_8(p_2)
&
y_4(p_3)\dot s_5y_6(p_5)x_7(m_6)\dot s_7^{-1} \dot s_2y_3(p_8)\dot s_4y_5(p_{10})y_6(p_{11})\\
&\cdot  
y_1(p_{12})x_2(m_{13})\dot s_2^{-1}y_3(p_{14})x_4(m_{15})\dot s_4^{-1}x_5(m_{16})\dot s_5^{-1}.
\end{align*}
The Go-diagram and the labeled Go-diagram are as follows:
\medskip
\[
\young[5,5,4,2][10][\hskip0.4cm\circle*{5},\hskip0.4cm\circle*{5},,\hskip0.4cm\circle*{5},,,,\hskip0.4cm\circle{5},,\hskip0.4cm\circle{5},\hskip0.4cm\circle*{5},,\hskip0.4cm\circle{5},,
,\hskip0.4cm\circle{5}]\hskip2cm \young[5,5,4,2][10][$-1$,$-1$,$p_{14}$,$-1$,$p_{12}$,
$p_{11}$,$p_{10}$,$1$,$p_8$,$1$,$-1$,$p_5$,$1$,$p_3$,$p_2$,$1$]
\]
\medskip\noindent
The $A$-matrix is then given by
\[
A=\begin{pmatrix}
-p_{12}p_{14}& q_{13} & p_{14} & q_{15} & -m_{16} & 1 & 0 & 0 & 0 \\
0 & p_8p_{10}p_{11} & 0 & p_{11}(p_3+p_{10}) & p_{11} & 0 & 1 & 0 & 0 \\
0  &  0 &  0  & -p_3p_5 & -p_5 & 0 & -m_6 & 1 & 0 \\
0 & 0 & 0 & 0 & 0 & 0 & p_2 & 0 & 1
\end{pmatrix},
\]
where the matrix entry $q_{13}=-m_{13}p_{14}+m_{15}p_8-m_{16}p_8p_{10}$ and $q_{15}=m_{15}-m_{16}(p_3+p_{10})$.  
In Figure \ref{fig:CP}, we show
contour plots $\CC_t(u_A)$ for the solution $u_A(x,y,t)$ at $t=-10, 0, 10$,
using the choice of 
parameters 
$(\kappa_1,\ldots,\kappa_9)=(-5,-3,-2,-1,0,1,2,3,4)$,  
$p_j=1$ for all $j$, and $m_l=0$ for all $\ell$. Note that:

\begin{itemize}
\item[(a)] For $y\gg 0$, there are four unbounded 
line-solitons, whose types from right to left are:
\[
[1,6],\quad [2,7],\quad [4,8],\quad {\rm and}\quad [7,9]
\]
\item[(b)] For $y\ll 0$, there are five unbounded 
line-solitons, whose types  from left to right are:
\[
[1,3],\quad [2,5],\quad [3,6],\quad [4,8],\quad {\rm and}\quad [5,9]
\]
\end{itemize}
We can see from this example  that the line-solitons for $y\gg0$ correspond to the excedances 
in $\pi=
(6,7,1,8,2,3,9,4,5)$, while 
those for $y\ll0$ correspond to the nonexcedances.
\begin{figure}[h]
\begin{center}
\includegraphics[height=1.9in]{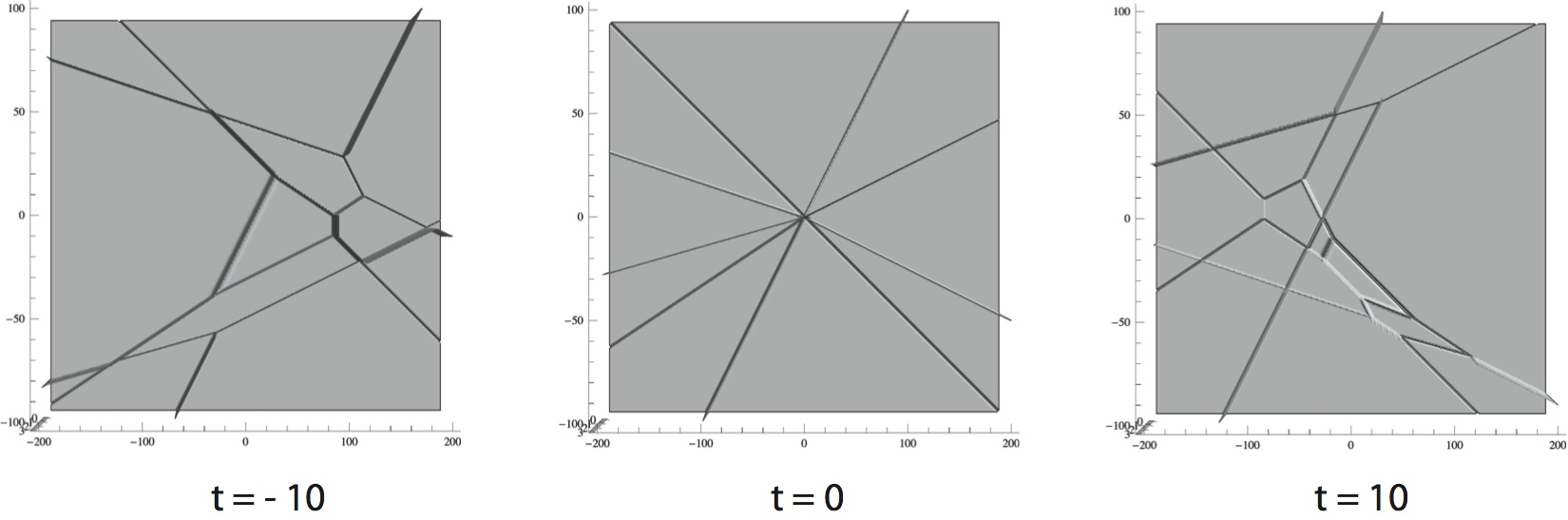}
\end{center}
\caption{Example of contour plots $\CC_t(u_A)$ for
$A\in Gr_{4,9}$.
The contour plots are obtained by ``Plot3D'' of Mathematica 
(see the details in the text).
 \label{fig:CP}}
\end{figure}

Note that if there are two adjacent regions of the contour plot
whose Pl\"ucker coordinates have different signs, then the line-soliton
separating them is singular.  
For example,
the line-soliton of type $[4,8]$ (the second soliton from the left in $y\gg0$) is 
singular, because the Pl\"ucker coordinates corresponding 
to the (dominant exponentials of the) adjacent regions are 
\[
\Delta_{1,2,4,9}=p_3p_5p_8p_{10}p_{11}p_{12}p_{14}=1\qquad {\rm and}\qquad
\Delta_{1,2,8,9}=-p_8p_{10}p_{11}p_{12}p_{14}=-1.
\]
\end{example}


\section{Soliton graphs and generalized plabic graphs}\label{sec:solitonplabic}

The following notion of {\it soliton graph}
forgets the metric
data of the contour plot, but preserves
the data of how line-solitons interact and which exponentials dominate.

\begin{definition}\label{soliton-graph}
Let $A \in Gr_{k,n}$ and consider a generic contour plot $\CC_{t}(u_A)$ 
for some time $t$.
Color a trivalent
vertex black (respectively, white)
if it has a unique edge extending downwards (respectively, upwards) from it.
We preserve the labeling of regions and edges that was used 
in the contour plot: we label a region by $E_I$ if the dominant 
exponential in that region is $\Delta_I E_I$, and label
each line-soliton by its {\it type} $[i,j]$ 
(see Lemma \ref{separating}).
We also preserve the topology of the graph,
but forget the metric structure.
We call this labeled graph with bicolored vertices
the \emph{soliton graph} $G_{t_0}(u_A)$.
\end{definition}

\begin{example}
We continue 
Example \ref{ex:Gr49CP}.
Figure \ref{contour-soliton} contains
the same contour plot $\CC_t(u_A)$ as that at the left of 
Figure \ref{fig:CP}.
One may use Lemma \ref{separating} to label all regions and edges
in the soliton graph.
After computing the Pl\"ucker coordinates, 
one can identify the singular solitons, which are indicated
by the dotted lines in the soliton graph.
\begin{figure}[h]
\begin{center}
\includegraphics[height=1.8in]{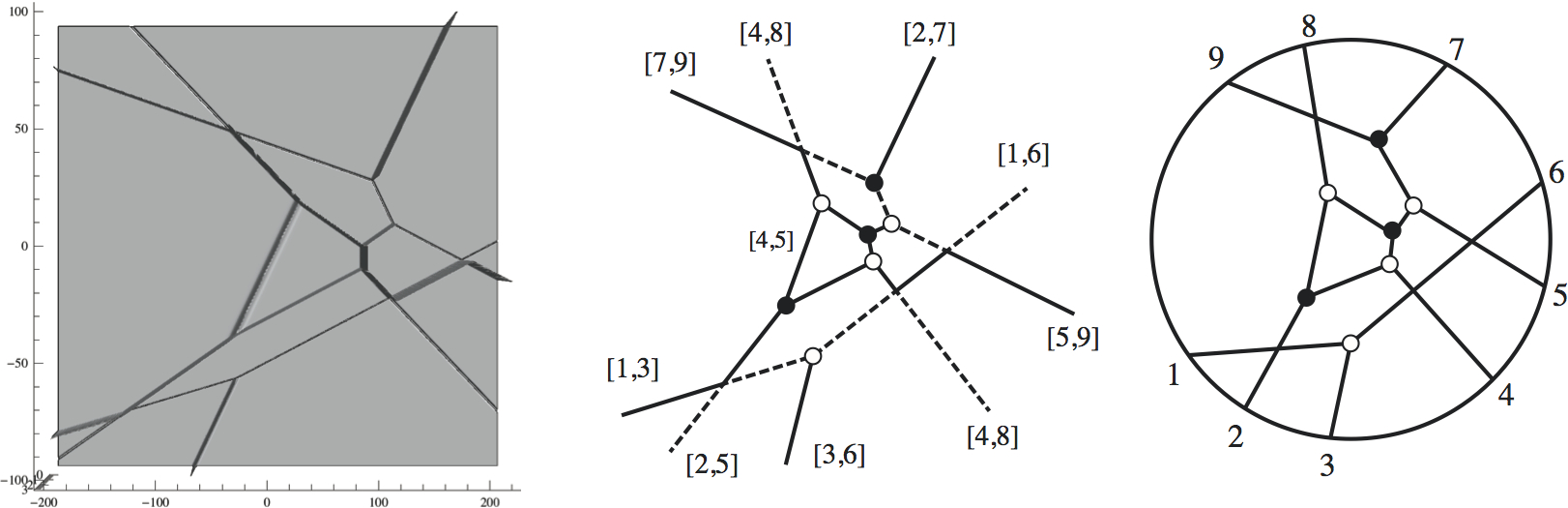}
\end{center}
\caption{Example of a contour plot $\CC_t(u_A)$, its
soliton graph $C=G_t(u_A)$, and its generalized plabic
graph $Pl(C)$.   The parameters used are
those from Example \ref{ex:Gr49CP}.
In particular,
$(\kappa_1,\dots,\kappa_9)=(-5,-3,-2,-1,0,1,2,3,4)$, and 
$\pi = (6,7,1,8,2,3,9,4,5).$
 \label{contour-soliton}}
\end{figure}
\end{example}


We now describe how to pass 
from a soliton graph to a \emph{generalized plabic graph}. 

\begin{definition}
A \emph{generalized plabic graph\/}
is an undirected graph $G$ drawn inside a disk
with $n$ \emph{boundary vertices\/} labeled
$\{1,\dots,n\}$.
We require that each boundary vertex $i$ is either isolated 
(in which case it is colored with color $1$ or $-1$), or
is incident
to a single edge; and each internal vertex is colored black or white.
Edges are allowed to cross each 
other in an $X$-crossing  (which is not considered to be a vertex). 
\end{definition}

By Theorem \ref{thm:soliton-perm}, the following construction
is well-defined.
\begin{definition}\label{soliton2plabic}
Fix a positroid stratum $\S_{\pi^{:}}$ of $Gr_{k,n}$ where
$\pi^{:} = (\pi,col)$.
To each soliton graph $C$ coming from a point of that stratum we associate 
a generalized plabic graph
$Pl(C)$ by:
\begin{itemize}
\item embedding $C$ into a disk, so that each unbounded line-soliton
of $C$ ends at a \emph{boundary vertex};
\item labeling the boundary vertex 
incident to the edge with labels $i$ and $\pi(i)$ 
by $\pi(i)$;
\item adding an isolated boundary vertex labeled $h$ with color
$1$ (respectively, $-1$) whenever $h\in J$ for each region label $E_J$
(respectively, whenever $h\notin J$ for any region label $E_J$);
\item forgetting the labels of all edges
and regions.
\end{itemize}
\end{definition}
See Figure \ref{contour-soliton} for a soliton graph $C$ together with the 
corresponding generalized plabic graph $Pl(C)$.

\begin{definition}\label{gen:trip}
Given a generalized plabic graph $G$,
the \emph{trip} $T_i$ is the directed path which starts at the boundary vertex 
$i$, and follows the ``rules of the road": it turns right at a 
black vertex,  left at a white vertex, and goes straight through the 
$X$-crossings.  
Note that $T_i$  will also 
end at a boundary vertex.  
If $i$ is an isolated vertex, then 
$T_i$ starts and ends at $i$.  
Define $\pi_G(i)=j$ whenever $T_i$
ends at $j$. 
It is not hard to show that $\pi_G$ is a permutation,
which we call the \emph{trip permutation}.
\end{definition}

We use the trips to label the edges and regions
of each generalized plabic graph.

\begin{definition}\label{labels}
Given a generalized plabic graph $G$, 
start at each non-isolated boundary
vertex $i$ and label every edge along trip $T_i$ with $i$.
Such a trip divides the disk containing $G$ into two parts: 
the part to the left of $T_i$, and the part to the right.
Place an $i$ in every region which is to the left of $T_i$.
If $h$ is an isolated boundary vertex with color $1$,
put an $h$ in every region of $G$.
After repeating this procedure for each boundary vertex,
each edge will be labeled by up to two numbers (between $1$ and $n$),
and each region will be labeled by a collection of numbers.
Two regions separated by an edge labeled by both $i$ and $j$ will have 
region labels $S$ and $(S\setminus \{i\}) \cup \{j\}$.
When an edge is labeled by two numbers $i<j$, we write $[i,j]$
on that edge, or $\{i,j\}$ or $\{j,i\}$ if we do not wish to specify
the order of $i$ and $j$.  
\end{definition}

\begin{theorem}\label{soliton-plabic}\cite[Theorem 7.6]{KW2}
Consider a soliton graph $C=G_t(u_A)$ coming from a point $A$
of a positroid stratum
$\S_{\pi^:}$, where $\pi^: = (\pi,col)$.
Then the trip permutation of $Pl(C)$ is $\pi$, 
and by labeling edges of $Pl(C)$ according
to Definition \ref{labels}, we will recover the original edge 
and region
labels in $C$.  
\end{theorem}

We invite the reader to verify Theorem \ref{soliton-plabic} for 
the graphs in Figure \ref{contour-soliton}.

\begin{remark} By Theorem \ref{soliton-plabic}, we can identify 
each soliton graph $C$ with its generalized plabic graph $Pl(C)$.
From now on, we will often ignore the labels of edges and regions
of a soliton graph, and simply record the labels on boundary vertices.
\end{remark}

\section{The contour plot  for $t\ll0$}\label{sec:t<<0}

Consider a matroid stratum $S_{\mathcal{M}}$ contained in 
the Deodhar component $S_{{D}}$,
where $D$ is the corresponding or Go-diagram.
From Definition \ref{contour} it is clear that the contour plot associated to 
any $A\in S_{\mathcal{M}}$ depends only on $\mathcal{M}$, not on $A$.  
In fact for $t\ll0$ a stronger statement is true -- 
the contour plot for 
any $A \in S_{\mathcal{M}} \subset S_D$ depends only on $D$, 
and not on $\mathcal{M}$.
In this section
we will explain how to use $D$ to construct first a 
generalized plabic graph $G_-(D)$,
and then the contour plot $\CC_{t}(\mathcal{M})$ for $t\ll 0.$

\subsection{Definition of the contour plot for $t \ll0$.}\label{sec:contour}
Recall from \eqref{f-contour} the definition of $f_{\mathcal{M}}(x,y,t)$.
To understand how it behaves for $t \ll 0$, let us
rescale everything by $t$.
Define $\bar{x} = \frac{x}{t}$ and $\bar{y} = \frac{y}{t}$,
and set
\[
\phi_i(\bar{x}, \bar{y}) = \kappa_i \bar{x} + 
\kappa_i^2 \bar{y} + \kappa_i^3,
\]
that is, $\kappa_i x + \kappa_i^2 y + \kappa_i^3 t = t\phi_i(\bar{x},\bar{y})$.
Note that because $t$ is negative,  $x$ and $y$ have the opposite
signs of $\bar{x}$ and $\bar{y}$.  This
leads to the following definition of the contour plot
for $t \ll 0$.

\begin{definition}\label{contourplot-infinity}
We define the contour plot
$\CC_{-\infty}(\mathcal{M})$
to be the locus in $\R^2$ where
\begin{equation*}
\underset{J\in\mathcal{M}}\min \left\{
                     \sum_{i=1}^k \phi_{j_i}(\bar{x},\bar{y}) \right\}\quad 
\text{ is not linear .}
\end{equation*}
\end{definition}

\begin{remark}\label{rem:rotate}
After a $180^{\circ}$ rotation,  $\CC_{-\infty}(\mathcal{M})$ is the limit of
 $\CC_{t}(u_A)$  as $t\to -\infty$, for any $A\in S_{\mathcal{M}}$.
Note that the rotation is required because the positive $x$-axis
(respectively, $y$-axis)
corresponds to the negative $\bar{x}$-axis (respectively,
$\bar{y}$-axis).
\end{remark}

\begin{definition}\label{def:v}
Define $v_{i,\ell,m}$ to be the point in $\R^2$ where
$\phi_i(\bar{x},\bar{y}) = 
\phi_{\ell}(\bar{x},\bar{y}) = 
\phi_m(\bar{x},\bar{y}).$
A simple calculation yields that 
the point $v_{i,\ell,m}$ has the following
coordinates in the $\bar{x}\bar{y}$-plane:
\[
v_{i,\ell,m}=(\kappa_i \kappa_{\ell} + \kappa_i \kappa_m +
\kappa_{\ell} \kappa_m, -(\kappa_i+\kappa_{\ell}+\kappa_m)).
\]
\end{definition}
Some of the points $v_{i,\ell,m}\in \mathbb{R}^2$ correspond to trivalent vertices in the contour plots we construct; such a point is the location
of the resonant interaction of three line-solitons of types $[i,\ell]$, $[\ell,m]$ and $[i,m]$
(see Theorem \ref{t<<0} below).
Because of our assumption on the genericity of the $\kappa$-parameters,
those points are all distinct.

\subsection{Main results on the contour plot for $t \ll0 $}

\begin{figure}[h]
\centering
\includegraphics[height=3.3in]{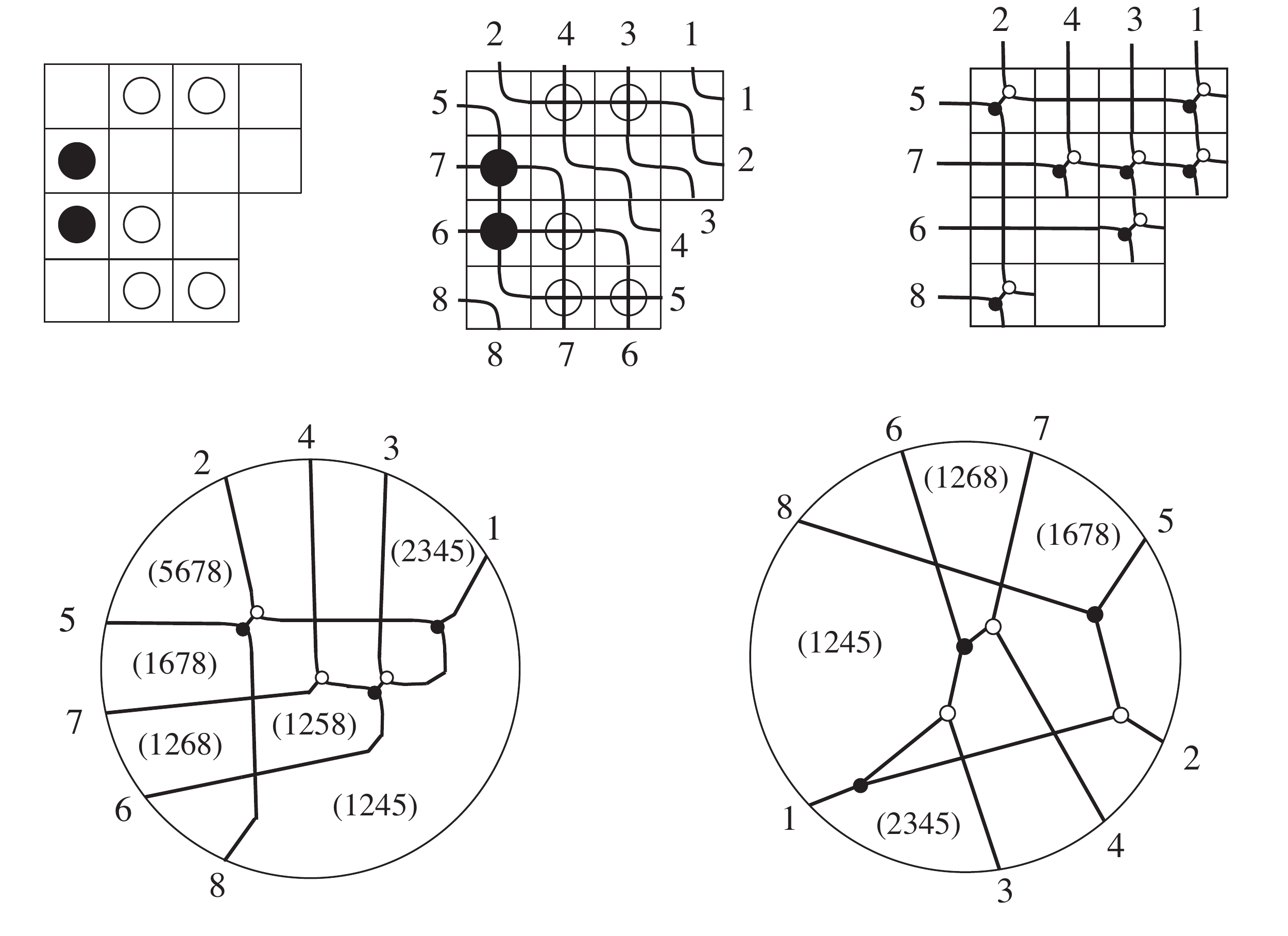}
\caption{Construction of the generalized plabic graph $G_-(D)$ associated to the Go-diagram $D$.  The labels of the regions of the graph indicate the index sets
of the corresponding Pl\"ucker coordinates.
Using the notation of Definition \ref{def:DtoPi}, we have
$\pi(D)=vw^{-1}=(5,7,1,6,8,3,4,2)$.}
\label{GoPlabic}
\end{figure}

This section is based on \cite[Section 10]{KW3},
and generalizes the results of  
\cite[Section 8]{KW2} to a soliton solution coming from an arbitrary point of the real Grassmannian (not just the non-negative part).
We start by giving an algorithm to construct a 
generalized plabic graph $G_-(D)$,
which will be used to construct  $\CC_{-\infty}(\mathcal{M})$.
Figure \ref{GoPlabic} illustrates the steps of Algorithm \ref{GoToPlabic},
starting
from the Go-diagram of the Deodhar component
$S_{D}$ where $D$ is as in the upper left corner of Figure \ref{GoPlabic}.

\begin{algorithm}  \label{GoToPlabic}
\cite[Algorithm 10.4]{KW3}
 From a Go-diagram $D$ to $G_-(D)$:
\begin{enumerate}
\item Start with a Go-diagram $D$ contained in a $k\times (n-k)$
rectangle, and
replace each \wstn, \bstn, and blank box by a cross, a cross, and a
pair of {\it elbows},
respectively.  Label the $n$ edges along the southeast border 
of the Young diagram by the numbers $1$ to $n$,
from northeast to southwest.  The configuration of crosses and elbows forms $n$ ``pipes" 
which travel from the southeast border to the northwest border; label the endpoint of 
each pipe by the label of its starting point.
\item Add a pair of black and white vertices to each pair of elbows,
and connect them by an edge,
as shown
in the upper right of Figure \ref{GoPlabic}.  Forget the labels
of the southeast border.  If there is an endpoint of a pipe on the east or south border whose pipe
starts by going straight, then erase the straight portion preceding the first elbow.
If there is a horizontal (respectively, vertical) pipe starting at $i$ with no elbows,
then erase it, and add an isolated boundary vertex labeled $i$
with color $1$ (respectively, $-1$).
\item Forget any degree $2$ vertices, and forget
any edges of the graph which end
at the southeast border of the diagram.
Denote the resulting
graph  $G_-(D)$.
\item After embedding the graph in a disk
with $n$ boundary vertices (including isolated vertices)
we obtain a generalized plabic graph,
which we also denote $G_-(D)$.
If desired, stretch and rotate $G_-(D)$ so that the boundary vertices
at the west side of the diagram are at the north instead. 
\end{enumerate}
\end{algorithm}

The following is the main result of this section.   
\begin{theorem}\label{t<<0}\cite[Theorem 10.6]{KW3}
Choose a matroid stratum $S_{\mathcal{M}}$ and let $S_D$ be the 
Deodhar component containing $S_{\mathcal{M}}$.
Recall the definition of $\pi(D)$ from 
Definition \ref{def:DtoPi}.
Use Algorithm \ref{GoToPlabic} 
to obtain $G_-(D)$.  
Then $G_-(D)$ has trip permutation $\pi(D)$, and we can use it to explicitly
construct $\CC_{-\infty}(\M)$ as follows.
Label the edges of $G_-(D)$
according to the rules of the road.  Label by $v_{i,\ell,m}$
each trivalent
vertex which is incident to edges labeled $[i,\ell]$, $[i,m]$, and $[\ell,m]$,
and give that vertex the coordinates
$(\bar{x},\bar{y}) = (\kappa_i \kappa_{\ell}+\kappa_i \kappa_m + \kappa_{\ell} \kappa_m,
-(\kappa_i+\kappa_{\ell}+\kappa_m))$.
Replace each edge labeled $[i,j]$ which ends 
at a boundary vertex by an unbounded line-soliton 
with slope $\kappa_i + \kappa_j$.
(Each edge labeled $[i,j]$ between two trivalent vertices 
will automatically
have slope
$\kappa_i + \kappa_j$.)
In particular, $\CC_{-\infty}(\M)$ is determined by $D$.
Recall from Remark \ref{rem:rotate} that after a
$180^{\circ}$ rotation,  $\CC_{-\infty}(\M)$ is the limit
of $\CC_{t}(u_A)$ as $t\to -\infty$, for any $A\in S_{\mathcal{M}}$.
\end{theorem}
\begin{remark}
Since the contour plot $\CC_{-\infty}(\M)$ depends only on $D$, we also 
refer to it as $\CC_{-\infty}(D)$.
\end{remark}

\begin{remark}\label{rem:t>>0}
The results of this section may be extended to the case $t \gg0$
by duality considerations (similar to the way in which
our previous paper \cite{KW2} described 
contour plots for both $t\ll0$ and $t\gg0$).  Note that 
the Deodhar decomposition of $Gr_{k,n}$ depends on a choice
of ordered basis $(e_1,\dots,e_n)$.  Using the ordered basis
$(e_n,\dots,e_1)$ instead and the corresponding
Deodhar decomposition, one may explicitly describe
contour plots at $t\gg0$.
\end{remark}

\begin{remark}\label{rem:differ}
Depending on the choice of the parameters $\kappa_i$,
the contour plot 
$\CC_{-\infty}(D)$ may have a slightly different topological structure than the soliton
graph $G_-(D)$.  While the incidences of line-solitons
with trivalent vertices are determined by $G_-(D)$, the locations 
of $X$-crossings may vary based on the $\kappa_i$'s. More specifically,
changing the $\kappa_i$'s may change the contour plot via a sequence
of \emph{slides}, see \cite[Section 11]{KW3}.
\end{remark}

\subsection{$X$-crossings in the contour plots}

Recall the notions of black and white $X$-crossings from 
Definition \ref{def:blackwhiteX}.
In \cite[Theorem 9.1]{KW2}, we proved that the presence of 
$X$-crossings in contour plots at $|t|\gg0$ implies 
that there is a two-term Pl\"ucker relation.

\begin{theorem}\cite[Theorem 9.1]{KW2}\label{2term}
Suppose that there is an 
$X$-crossing in a contour plot $\CC_t(u_A)$ 
for some $A \in Gr_{k,n}$ where $|t|\gg0$.  
Let $I_1$, $I_2$, $I_3$, and $I_4$ be the $k$-element subsets of $\{1,\dots,n\}$
corresponding to the dominant exponentials incident
to the $X$-crossing listed in circular order.  
\begin{itemize}
\item If the $X$-crossing is white, 
we have $\Delta_{I_1}(A) \Delta_{I_3}(A) = \Delta_{I_2}(A) \Delta_{I_4}(A).$
\item If
the $X$-crossing is black, 
we have $\Delta_{I_1}(A) \Delta_{I_3}(A) = -\Delta_{I_2}(A) \Delta_{I_4}(A).$
\end{itemize}
\end{theorem}

The following corollary is immediate.
\begin{corollary}\label{prop:opposite}
If there is a black $X$-crossing in a contour plot at $t\ll 0$ or $t\gg 0$, 
then among the Pl\"ucker coordinates associated to 
the dominant exponentials incident 
to that black $X$-crossing, three must be positive and one negative,
or vice-versa. 
\end{corollary}

\begin{corollary}\label{cor:white}
Let $D$ be a $\Le$-diagram, that is, a Go-diagram with no black stones.
Let $A \in S_D$ and $t \ll0$.  Choose any $\kappa_1 < \dots < \kappa_n$.
Then the contour plot $\CC_t(u_A)$ can have only white $X$-crossings.
\end{corollary}

\section{Total positivity, regularity, and cluster algebras}\label{sec:regularity}

In this paper we have been studying 
solutions $u_A(x,y,t)$ to the KP equation coming from points of the real Grassmannian.
Among these solutions, those coming from the totally non-negative part of 
$Gr_{k,n}$ are especially nice.  In particular, $(Gr_{k,n})_{\geq 0}$
parameterizes the \emph{regular} soliton solutions coming from $Gr_{k,n}$,
see Theorem \ref{th:regularity}.  
This result provides an important motivation for studying
the soliton solutions coming from 
$(Gr_{k,n})_{\geq 0}$.  
After discussing the regularity result below, we will discuss \emph{positivity tests}
for elements of the Grassmannian, as well as the connection between 
soliton solutions from $(Gr_{k,n})_{>0}$ and the theory of cluster algebras.

\begin{theorem}\label{th:regularity}\cite[Theorem 12.1]{KW3}
Fix parameters $\kappa_1 < \dots < \kappa_n$ and 
an element $A \in Gr_{k,n}$.  Consider the corresponding soliton solution
$u_A(x,y,t)$ of the KP equation.  This solution is regular at $t\ll0$
if and only if $A \in (Gr_{k,n})_{\geq 0}$.
Therefore this solution is regular for all times $t$ 
if and only if $A \in (Gr_{k,n})_{\geq 0}$.
\end{theorem}

If one is studying total positivity on the Grassmannian, then a natural 
question is the following: given $A\in Gr_{k,n}$, how many Pl\"ucker 
coordinates, and which ones, must one compute, in order to determine
that $A \in (Gr_{k,n})_{\geq 0}$?
This leads to the following notion of \emph{positivity test}.

\begin{definition}\label{def:pos-test}
Consider the Deodhar component $S_D \subset Gr_{k,n}$, where $D$
is a Go-diagram.
A collection $\J$ of $k$-element subsets of $\{1,2,\dots,n\}$
is called a \emph{positivity test} for $S_D$ if for any $A \in S_D$,
the condition that $\Delta_I(A) >0$ for all $I \in \J$ 
implies that $A \in (Gr_{k,n})_{\geq 0}$.
\end{definition}

It turns out that the collection of Pl\"ucker coordinates
corresponding to dominant exponentials in contour plots
at $t\ll0$ provide positivity tests for positroid strata.

\begin{theorem}\label{th:Le}\cite[Theorem 12.9]{KW3}
Let $A \in S_D \subset Gr_{k,n}$, where $D$ is a $\Le$-diagram, and let
$t \ll0$.  If $\Delta_J(A)>0$ for each dominant exponential $E_J$
in the contour plot $\CC_t(u_A)$,
then  $A \in (Gr_{k,n})_{\geq 0}$.  In other words,
the Pl\"ucker coordinates corresponding to the 
dominant exponentials in $\CC_t(u_A)$ comprise a positivity test for 
$S_D$.
\end{theorem}

\begin{remark} Recall that the Deodhar components 
$S_D$ have non-empty intersection with $(Gr_{k,n})_{\geq 0}$
unless $D$ is a $\Le$-diagram.  Therefore Theorem \ref{th:Le}
restricts to the case that $D$ is a $\Le$-diagram.
\end{remark}


\subsection{TP Schubert cells and reduced plabic graphs}\label{Reduced-Cluster}

In this section we provide a new characterization of the so-called
\emph{reduced} plabic graphs \cite[Section 12]{Postnikov}.  This will
allow us to make a connection to cluster algebras in Section \ref{sec:cluster}.

We will always assume that a plabic graph is {\it leafless}, i.e. that 
it has no non-boundary leaves, and that it has no isolated components.
In order to define {\it reduced}, we first
define some local transformations of plabic graphs.

(M1) SQUARE MOVE.  If a plabic graph has a square formed by
four trivalent vertices whose colors alternate,
then we can switch the
colors of these four vertices.

(M2) UNICOLORED EDGE CONTRACTION/UNCONTRACTION.
If a plabic graph contains an edge with two vertices of the same color,
then we can contract this edge into a single vertex with the same color.
We can also uncontract a vertex into an edge with vertices of the same
color.

(M3) MIDDLE VERTEX INSERTION/REMOVAL.
If a plabic graph contains a vertex of degree 2,
then we can remove this vertex and glue the incident
edges together; on the other hand, we can always
insert a vertex (of any color) in the middle of any edge.

\begin{figure}[h]
\centering
\includegraphics[height=.5in]{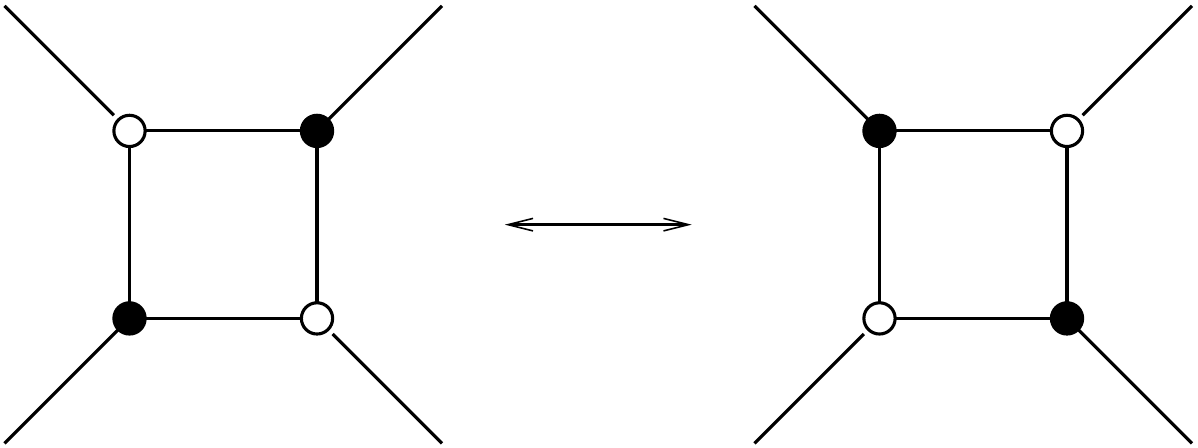}
\hspace{1cm}
\raisebox{0.2cm}{\includegraphics[height=.3in]{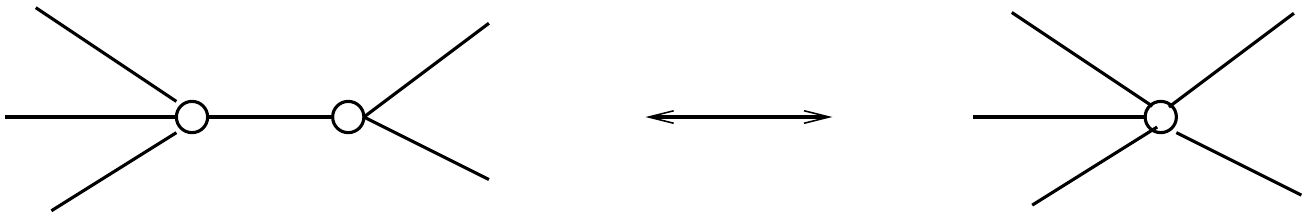}}
\hspace{1cm}
\raisebox{0.5cm}{\includegraphics[height=.06in]{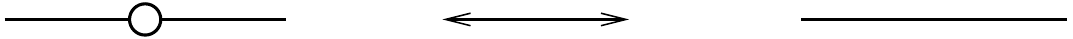}}
\caption{A square move; a unicolored edge contraction;  a middle vertex insertion/ removal}
\label{M}
\end{figure}

(R1) PARALLEL EDGE REDUCTION.  If a network contains
two trivalent vertices of different colors connected
by a pair of parallel edges, then we can remove these
vertices and edges, and glue the remaining pair of edges together.

\begin{figure}[h]
\centering
\includegraphics[height=.25in]{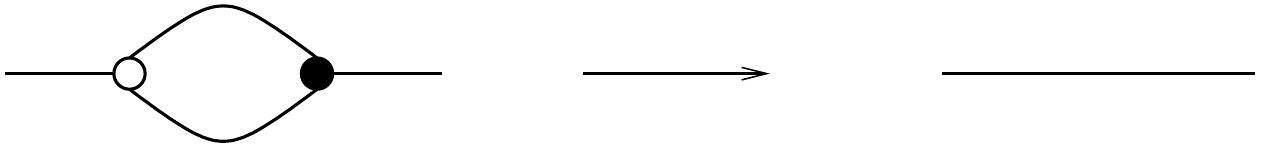}
\caption{Parallel edge reduction}
\label{R1}
\end{figure}

\begin{definition}\cite{Postnikov}
Two plabic graphs are called \emph{move-equivalent} if they can be obtained
from each other by moves (M1)-(M3).  The \emph{move-equivalence class} 
of a given plabic graph $G$ is the set of all plabic graphs which are move-equivalent
to $G$.
A leafless plabic graph without isolated components
is called \emph{reduced} if there is no graph in its move-equivalence 
class to which we can apply (R1).
\end{definition}

\begin{theorem}\cite[Theorem 13.4]{Postnikov}
Two reduced plabic graphs which each have $n$ boundary vertices
are move-equivalent if and only if they have the same 
trip permutation.
\end{theorem}

Our new characterization of reduced plabic graphs is as follows.

\begin{definition}\label{def:resonance}
We say that a (generalized) plabic graph has the 
\emph{resonance property}, if after labeling edges via Definition 
\ref{labels}, the set $E$ of edges incident to a given vertex  has
the following property:
\begin{itemize}
\item  there exist numbers $i_1<i_2<\dots<i_m$ such that when 
we read the labels of $E$,  we see the labels
$[i_1,i_2],[i_2,i_3],\dots,[i_{m-1},i_m],[i_1,i_m]$ appear in 
counterclockwise order.
\end{itemize}
\end{definition}

\begin{theorem}\label{th:reduced}\cite[Theorem 10.5]{KW2}.
A plabic graph is reduced 
if and only if it has the resonance property.\footnote{Recall 
from Definition \ref{def:plabic} that our convention is to label 
boundary vertices of a plabic graph $1,2,\dots,n$ in counterclockwise
order.  If one chooses the opposite convention,
then one must replace the word \emph{counterclockwise} in 
Definition \ref{def:resonance}
by \emph{clockwise}.}
\end{theorem}

\begin{corollary}\label{cor:reduced}\cite[Corollary 10.9]{KW2}
Suppose that $A$ lies in a TP Schubert cell, and that for some time 
$t$, $G_t(u_A)$ is generic with no X-crossings.  Then 
$G_t(u_A)$ is a reduced plabic graph.
\end{corollary}

\subsection{The connection to cluster algebras}\label{sec:cluster}
Cluster algebras are a class of commutative rings, introduced
by Fomin and Zelevinsky \cite{FZ}, which have a 
remarkable combinatorial structure.
Many coordinate rings of homogeneous spaces have
a cluster algebra structure: as shown by Scott \cite{Scott},
the Grassmannian is one such example.

\begin{theorem} \cite{Scott}\label{cluster-Grassmannian}
The coordinate ring of 
(the affine cone over) $Gr_{k,n}$  has  a 
cluster algebra structure.  Moreover, the set of Plucker coordinates
whose indices come from the labels of the regions
of a reduced plabic graph for $(Gr_{k,n})_{>0}$ 
comprises a \emph{cluster} for this cluster algebra.
\end{theorem}

\begin{remark}
In fact \cite{Scott} used the combinatorics of 
{\it alternating strand diagrams}, not reduced plabic graphs,
to describe clusters.
However, alternating strand 
diagrams are 
in bijection with reduced plabic graphs \cite{Postnikov}.
\end{remark}

\begin{theorem}\label{soliton-cluster}\cite[10.12]{KW2}
The set of Pl\"ucker coordinates labeling regions of a
generic trivalent soliton graph for the TP Grassmannian
is a cluster for the cluster algebra
associated to the Grassmannian.
\end{theorem}

Conjecturally, every positroid cell $\S_{\pi}^{tnn}$
of the totally non-negative
Grassmannian also carries a cluster algebra structure, 
and the Pl\"ucker coordinates labeling the regions of any reduced
plabic graph for $\S_{\pi}^{tnn}$ should be a cluster for that cluster
algebra.  In particular, the TP Schubert cells should carry cluster algebra
structures.  Therefore we conjecture that  Theorem
\ref{soliton-cluster} holds with ``Schubert cell" replacing 
``Grassmannian."  Finally, there should be a suitable generalization
of Theorem \ref{soliton-cluster} for arbitrary 
positroid cells.

\subsection{Soliton graphs for $(Gr_{2,n})_{>0}$ and triangulations
of an $n$-gon}

In this section we must use a slightly more general setup
for $\tau$-functions $\tau_A(t_1,\dots,t_m)$
and soliton solutions $u_A(t_1,\dots,t_m)$, 
in which the exponential functions
$E_i$ and $E_I$ are functions of variables $t_1,\dots,t_m$, 
and $m \geq n$.  Here $t_1=x$, $t_2=y$, $t_3=t$, and the 
other $t_i$'s are referred to as \emph{higher times}.
See \cite[Section 3]{KW2}
for details.
\begin{figure}[h]
\centering
\includegraphics[height=2.2in]{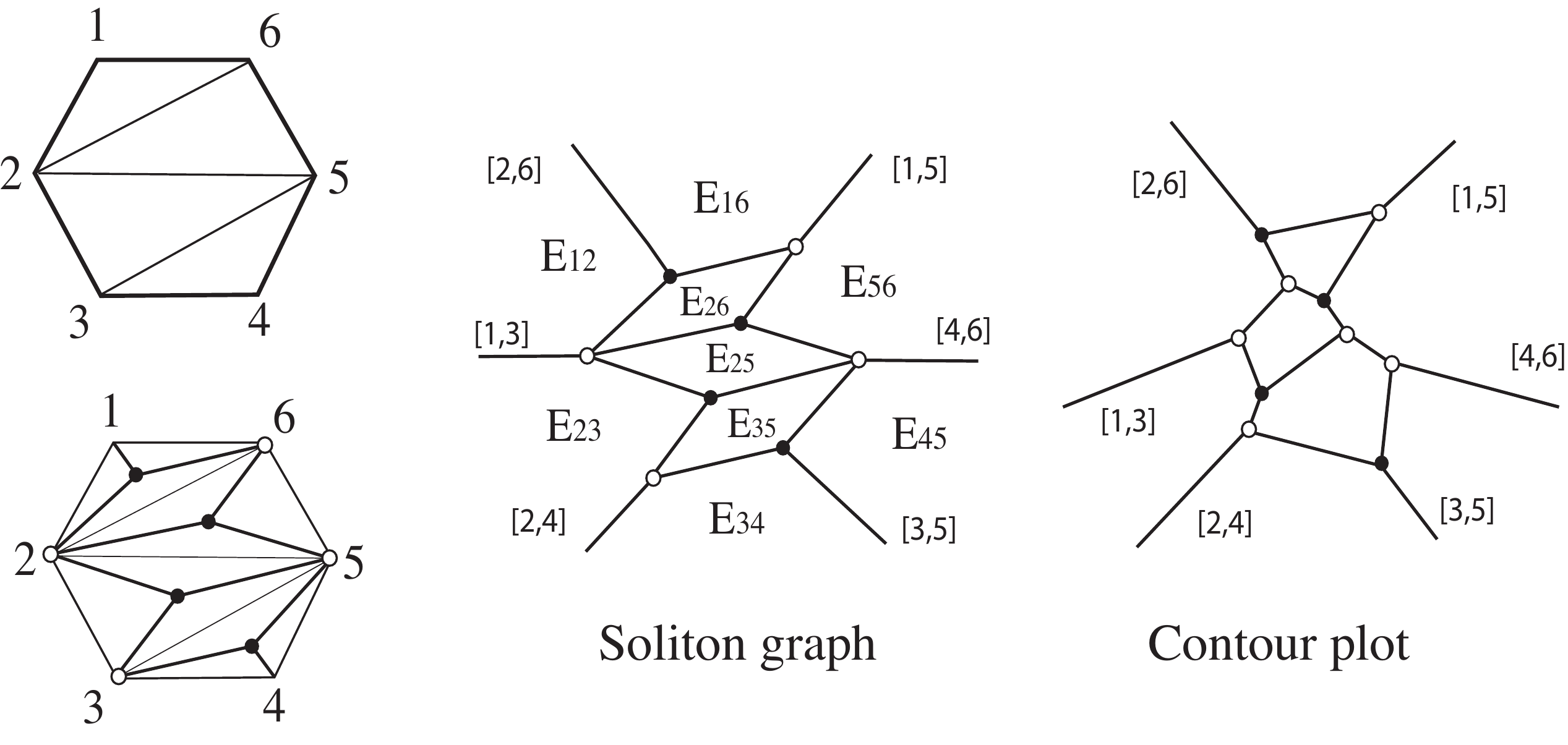}
\caption{Algorithm \ref{TriangulationA}, starting from a triangulation of
a hexagon.  The right figure shows the corresponding contour plot for the soliton solution
for  $(Gr_{2,6})_{>0}$.
\label{Psi}}
\end{figure}

\begin{algorithm}  \label{TriangulationA} \cite[Algorithm 12.1]{KW2}
Let $T$ be a triangulation of an $n$-gon $P$, whose $n$ vertices are labeled
by the numbers $1,2,\dots,n$, in counterclockwise order. Therefore
each edge of $P$ and each diagonal of $T$ is specified by a pair of distinct integers
between $1$ and $n$. The following
procedure yields a labeled graph  $\Psi(T)$.
\begin{enumerate}
\item Put a black vertex in the interior of each triangle in $T$.
\item Put a white vertex at each of the $n$ vertices of $P$
which is incident to a diagonal of $T$; put a black vertex
at the remaining vertices of $P$.
\item Connect each vertex which is
inside a triangle of $T$ to the three vertices
of that triangle.
\item Erase the edges of $T$, and contract every pair of  adjacent vertices
which have the same color.  This produces a new graph $G$
with $n$ boundary vertices, in bijection with the vertices of the original
$n$-gon $P$.
\item Add one unbounded ray to each of the boundary vertices of
$G$, so as to produce a new (planar) graph  $\Psi(T)$.  Note that
$\Psi(T)$ divides the plane into regions; the bounded regions correspond to
the diagonals of $T$, and the unbounded regions correspond to the edges of $P$.
\end{enumerate}
\end{algorithm}
\begin{figure}[h]
\centering
\includegraphics[height=3.5in]{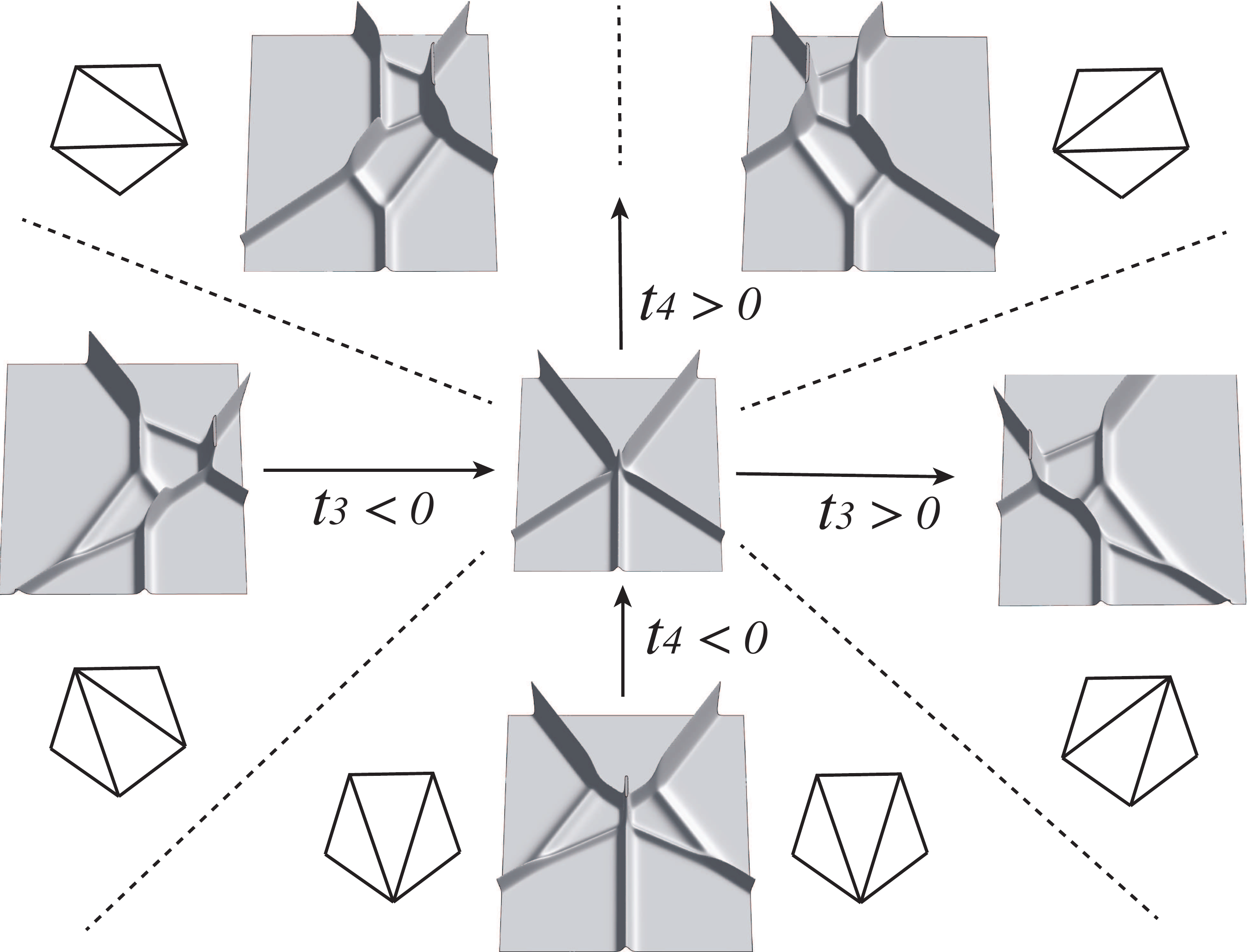}
\caption{Contour plots corresponding to the different soliton graphs for 
$(Gr_{2,5})_{>0}$. The center figure shows the contour plot at $t_3=t_4=0$.
\label{fig:G25}}
\end{figure}

See Figure \ref{Psi}.
Our main result in this section is the following.
\begin{theorem}\label{theorem:maintriangulation}\cite[Theorem 12.1]{KW2}
The graphs $\Psi(T)$ constructed above {\it are} soliton graphs
for $(Gr_{2,n})_{>0}$, and conversely, up to (M2)-equivalence,
any trivalent generic soliton graph
for $(Gr_{2,n})_{>0}$ comes from this construction.
Moreover, one can realize each graph $\Psi(T)$ by either:
\begin{itemize}
\item choosing an arbitrary $A \in (Gr_{2,n})_{>0}$ and varying the higher
times $t_3,\dots,t_n$ appropriately, or 
\item fixing an arbitrary collection of higher times $t_3,\dots,t_n$,
and using the torus action to choose an appropriate $A \in (Gr_{2,n})_{>0}$.
\end{itemize}
\end{theorem}

Figure \ref{fig:G25} shows the five triangulations of a 
pentagon, together with (contour plots corresponding 
to) the five different soliton graphs
which one may obtain from $(Gr_{2,5})_{>0}$.

\begin{remark}
The process of flipping a diagonal in the triangulation corresponds
to a mutation in the cluster algebra.  In the terminology of reduced plabic
graphs, a mutation corresponds to the square move (M1).
In the setting of KP solitons, each mutation may be considered as an
evolution along a particular flow of the KP hierarchy
defined by the symmetries of
the KP equation.
\end{remark}

\begin{remark}
It is known already that the set of reduced plabic graphs for the TP part of
$Gr_{2,n}$ all have the form given by Algorithm \ref{TriangulationA}.  And 
by Corollary \ref{cor:reduced}, every generic soliton graph is a reduced plabic
graph.  Therefore
it follows immediately that every
soliton graph for the TP part of $Gr_{2,n}$ must have the form of Algorithm \ref{TriangulationA}.
To prove 
Theorem \ref{theorem:maintriangulation}, one must also show that 
every outcome of Algorithm \ref{TriangulationA} can be realized as a soliton graph.
\end{remark}

\section{The inverse problem for soliton graphs}\label{sec:inverse}
\label{inverse}

The {\it inverse} problem for soliton solutions of the KP equation
 is the following:
given a time $t$ together with the contour plot 
$\CC_t(u_A)$ of a soliton 
solution, can one reconstruct the point $A$ of 
$Gr_{k,n}$ which gave rise to the solution?
Note that solving for $A$ is desirable, because this information would allow us to 
compute the entire past and the entire future of the soliton solution.


Using the cluster algebra structure 
for Grassmannians, we have the following.

\begin{theorem}\label{inverse2}\cite[Theorem 11.2]{KW2}
Consider a generic contour plot $\CC_t(u_A)$ of a soliton solution
which has no $X$-crossings, and which
comes from a point $A$ of the totally positive Grassmannian
at an {\it arbitrary} time $t$.
Then from the contour plot together with $t$ we 
can uniquely reconstruct the point $A$.
\end{theorem}

Using the description of contour plots of soliton solutions when 
$t \ll 0$, we have the following.

\begin{theorem}\label{inverse1}\cite[Theorem 11.3]{KW2}
Fix $\kappa_1<\dots<\kappa_n$ as usual.
Consider a generic contour plot $\CC_t(u_A)$ of a soliton solution
coming from a point $A$ of a positroid cell $S_{\pi}^{tnn}$,
for $t\ll0$.  Then from the contour plot together with $t$ we 
can uniquely reconstruct the point $A$.
\end{theorem}

\subsection{Non-uniqueness of the evolution of the contour plots for $t\gg 0$}
In contrast to the 
totally non-negative case,
where the soliton solution can be uniquely 
determined by the information in the contour plot at $t\ll 0$, 
if we consider arbitrary points $A\in Gr_{k,n}$, we cannot solve the 
inverse problem.

Consider $A \in S_D \subset Gr_{k,n}$.  
If the contour plot $\CC_{-\infty}(D)$
is topologically identical to $G_-(D)$, then 
the contour plot has almost no dependence 
on the parameters $m_j$ from the 
parameterization of $S_D$.  
This is because the Pl\"ucker coordinates
corresponding to the regions of 
$\CC_{-\infty}(D)$ (representing the dominant exponentials)
are either monomials in the $p_i$'s (see
\cite[Section 5]{KW2} and 
\cite[Remark 10.5]{KW2}),
or determined
from these by a ``two-term" Pl\"ucker relation.

Therefore it is possible to choose two different points $A$ and $A'$ in $S_D \subset Gr_{k,n}$
whose contour plots for a fixed $\kappa_1 < \dots \kappa_n$ and fixed $t \ll0$ are identical
(up to some exponentially small difference); 
we use the same parameters $p_i$
but different parameters $m_j$ for defining $A$ and $A'$.  However, as $t$ increases, those
contour plots may evolve to give different patterns.

\begin{figure}[h]
\begin{center}
\includegraphics[height=1.7in]{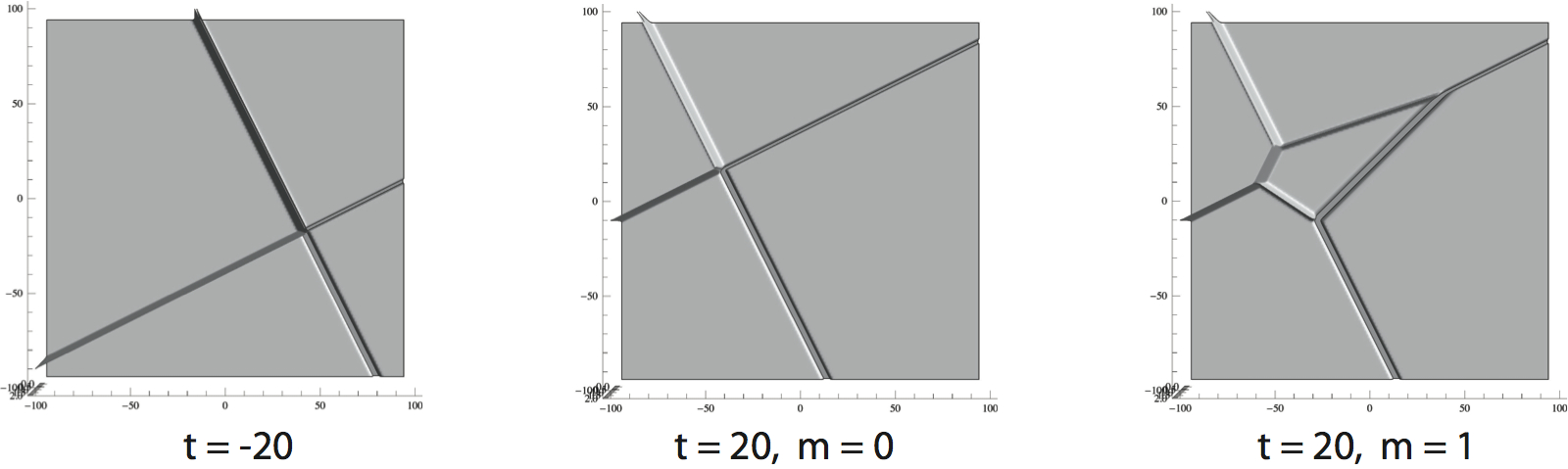}
\end{center}
\caption{The non-uniqueness of the evolution of the contour plots (and soliton graphs).  
The left panel shows the contour plot
at $t=-20$ for any value of $m$.  The middle panel shows the graph at $t=20$ with $m=0$,
and the right one shows the graph at $t=20$ with $m=1$.
These contour plots were made using the choice $p_i = 1$ for all $i$,
and 
$(\kappa_1,\ldots,\kappa_4)=(-2,-1,0,1.5)$.
In all of them,
the region at $x\gg 0$ has a positive sign ($\Delta_{3,4}=1$) and other
regions have negative signs. 
This means that the solitons adjacent to the region for $x\gg0$ are singular.}
 \label{fig:NUness}
\end{figure}

Consider the Deodhar stratum $S_D \subset Gr_{2,4}$, corresponding to
\[
{\bf w}=s_2s_3s_1s_2\text{ and }  {\bf v}=s_211s_2.
\]
The Go-diagram and labeled Go-diagram are given by
\[
\young[2,2][10][\hskip0.4cm\circle*{5},,,\hskip0.4cm\circle{5}]\qquad\qquad
\young[2,2][10][$-1$,$p_3$,$p_2$,$1$].
\]
\medskip
The matrix $g$ is calculated as $g=s_2y_3(p_2)y_1(p_3)x_2(m)s_2^{-1}$,
and its
projection to $Gr_{2,4}$ is 
\[
A=\begin{pmatrix}
-p_3 & -m & 1 & 0 \\
0 & p_2 & 0   & 1
\end{pmatrix}.
\]
The $\tau$-function is 
\[
\tau_A=-(p_2p_3E_{1,2}+p_3E_{1,4}+mE_{2,4}+p_2E_{2,3}-E_{3,4}),
\]
where $E_{i,j}:=(\kappa_j-\kappa_i)\exp(\theta_i+\theta_j)$.
The contour plots of the solutions with $m=0$ and $m\ne0$ are the same
(except for some exponentially small difference) 
when $t\ll0$.  In both cases, the plot consists of 
two line-solitons forming an $X$-crossing, where the parts of those solitons adjacent to the region
with dominant exponential $E_{3,4}$ (i.e. for $x\gg 0$) are singular,
see the left of Figure \ref{fig:NUness}.  

On the other hand, for $t\gg0$, 
the contour plot with $m=0$ is topologically the same as it was for $t\ll0$,
while the contour plot with $m\ne0$ has a box with 
dominant exponential $E_{2,4}$,
surrounded by four bounded solitons
(some of which are singular).  See the middle and right of 
Figure \ref{fig:NUness}.  
So not only the contour plots but also the soliton graphs are 
different for $t\gg0$!

Note that the non-uniqueness of the evolution of the 
contour plot (a tropical approximation) does not imply the non-uniqueness
of the evolution of the solution of the KP equation as $t$ changes.
If one makes two different choices for the $m_i$'s, 
the corresponding $\tau$-functions are different, but 
there is only an exponentially small difference in the corresponding
contour plots (hence the topology of the contour plots is identical).

\raggedright

\end{document}